\documentclass[preprint,10pt]{elsarticle}
\usepackage[utf8]{inputenc}
\usepackage{graphicx, float}
\usepackage{amsmath, amsfonts}
\usepackage{xcolor}
\usepackage{appendix}

\newcounter{todocounter}

\begin{document}

\begin{frontmatter}

\title{Differentiating densities on smooth manifolds}
\author{Adam A. \'Sliwiak\corref{cor}}
\ead{asliwiak@mit.edu}
\author{Qiqi Wang}
\ead{qiqi@mit.edu}
\address{Center for Computational Science and Engineering, Massachusetts Institute of Technology (MIT), 77 Massachusetts Avenue, Cambridge, MA, 02139, USA}
\address{MIT Department of Aeronautics and Astronautics}
\cortext[cor]{Corresponding author.}

\begin{abstract}
    Lebesgue integration of derivatives of strongly-oscillatory functions is a recurring challenge in computational science and engineering. Integration by parts is an effective remedy for huge computational costs associated with Monte Carlo integration schemes. In case of Lebesgue integrals over a smooth manifold, however, integration by parts gives rise to a derivative of the density implied by charts describing the domain manifold. This paper focuses on the computation of that derivative, which we call the {\em density gradient} function, on general smooth manifolds. We analytically derive formulas for the density gradient and present examples of manifolds determined by popular differential equation-driven systems. We highlight the significance of the density gradient by demonstrating a numerical example of Monte Carlo integration involving oscillatory integrands.  
\end{abstract}
\begin{keyword}
Density gradient function, Differentiable manifold, Lebesgue integral, Monte Carlo integration, Linear response 
\end{keyword}
\end{frontmatter}

%SECTION 1----------------------
\section{Introduction}
%Fundamental theorem of calculus
The fundamental theorem of calculus states that
$\int_a^b \partial_x f(x)\,dx = f(b) - f(a)$, where $f$ is some smooth function defined on the compact interval $[a,b]$. This theorem is critical in many applications, including computational sciences \cite{israel-theorem,thomson-theorem}. For example, if the function $f$ is strongly oscillatory, a numerical quadrature on the left-hand side would require many points and much computation to obtain accurate results.  Nevertheless, the fundamental theorem guarantees that the positive and negative derivatives of this oscillatory function largely cancel each other out.  Indeed, one can simply compute the right-hand side directly, without truncation error.

%Generalization of the fundamental theorem
As a generalization of the fundamental theorem, consider the integral of $\partial_x f$ over the same domain under some Lebesgue measure $m$, which is an antiderivative of the density function $\rho$ (a.k.a. the Radon-Nikodym derivative \cite{nagy-density}), i.e., $dm(x) = \rho(x)\,dx$. In the classical version of the theorem, as mentioned in the first paragraph, the density $\rho$ is constant and equals $1/(b-a)$ everywhere on the domain. If this is not the case, however, the integration by parts of $\partial_x f$ involves the derivative of $\rho$,  
\begin{equation}
\label{eqn:intro1}
\int_a^b \partial_x f(x)\,dm(x) = f\rho\Big|_a^b - \int_a^b f(x)\,\partial_x\rho(x)\,dx = f\rho\Big|_a^b - \int_a^b f(x)\,\frac{\partial_x\rho}{\rho}(x)\,dm(x)
\end{equation}
The integral in Eq. \ref{eqn:intro1} can be approximated using a Monte Carlo integration scheme if a set of realizations of $x$, $\{x^{1},x^{2},...,x^{N}\}$, distributed according to $m$, is given. However, if $f$ is a strongly-oscillatory function with large magnitude, the Monte Carlo method applied directly to the integral on the left-hand side (LHS) of Eq. \ref{eqn:intro1} would require a large amount of data to obtain an approximation with a reasonably small error \cite{olver-oscillatory,makri-oscillatory}. Alternatively, one can consider the right-hand side (RHS) of the same equation, which requires the function $f$ itself, not its derivative. Assuming the density $\rho$ is a well-behaved function, the variance of the integrand on the RHS is significantly smaller and, therefore, remarkably less data is needed to obtain an accurate result. However, extra computational effort must be put to evaluate $\partial_x\rho/\rho = \partial_x\log\rho$. The computation of that function, which we denote by $g$ and call it the {\em density gradient}, is the main focus of this paper. 

%Appearance of g in literature
Lebesgue integrals involving functions with high fluctuations are critical in the field of sensitivity analysis of chaotic dynamical systems. Ruelle \cite{ruelle-original,ruelle-corrections} derived a closed-form expression, known as the {\em linear response} formula, for the parametric derivative of the mean of a quantity of interest $J$. The linear response formula includes Lebesgue integrals of directional derivatives of a strongly oscillatory $J$ over the manifold of a chaotic system. A regularized version of Ruelle's formula, known as the space-split sensitivity (S3), was obtained through the integration by parts of the original formulation \cite{chandramoorthy-s3}. The S3 algorithm was successfully applied in various low-dimensional systems in the computation \cite{sliwiak-1d} and assessment of existence  \cite{sliwiak-differentiability} of parametric derivatives of statistical quantities describing chaos. The crux of the computation of the regularized Ruelle's formula is the {\em SRB density gradient}, defined as a directional derivative of the logarithm of the SRB density \cite{young-srb,crimmins-srb} along the unstable manifold. While an efficient numerical procedure for the approximation of the SRB density gradient specialized to systems with one-dimensional unstable manifolds is available \cite{chandramoorthy-s3,sliwiak-1d,sliwiak-differentiability, chandramoorthy-clv}, we still lack a generalizable algorithm applicable to arbitrary higher-dimensional chaotic systems.

%Focus of this paper
The main purpose of this work is to derive a general formula for the density gradient $g$, defined on a differentiable $m$-dimensional manifold $M$ immersed in the Euclidean space $\mathbb{R}^n$, $m\leq n$. In our analysis, we parameterize $M$ using the chart $x(\xi):\mathbb{R}^m\to \mathbb{R}^n$.  Here, the $g$ function is an $m$-element vector, where the $i$-th component equals a directional derivative of $\log\rho$, in the direction of a unit vector $s_{i}$, i.e., $g_{i} = \partial_{s_i}\rho/\rho = (\nabla_x\rho\cdot s_i)/\rho$. The scalar function $\rho$ is the density implied by $x(\xi)$. Without loss of generality, we assume that the $i$-th directional derivative is computed along the isoparametric line in the direction of increasing $i$-th component of $\xi$. Analogously to Eq. \ref{eqn:intro1}, the Lebesgue integral of the directional derivative of $J$ over $M$ with measure $m$ can be written using $g_i$,
\begin{equation}
    \label{eqn:intro2}
    \int_{M} \nabla_{x}J(x)\cdot s_{i}(x)\,dm(x) = - \int_{M} J(x)\,g_{i}(x)\,dm(x), 
\end{equation}
where $J$ is assumed to vanish on the boundary of $M$. For the reasons indicated above, it is computationally efficient to apply the Monte Carlo method to the RHS of Eq. \ref{eqn:intro2}. Analogous integration by parts is required to regularize the linear response \cite{chandramoorthy-s3}. Thus, the derivation of a computable expression for the density gradient defined on higher-dimensional smooth manifolds is a milestone in constructing algorithms for differentiating SRB measures. In addition, an explicit formula for $g$ might serve as a valuable tool in general numerical procedures involving integrals over geometrically complex domains. 

%Outline of the paper
The structure of this paper is the following. First, in Section \ref{sec:1D}, we derive a computable expression for the density gradient defined on one-dimensional manifolds (straight lines and curves). We also demonstrate a numerical example of Monte Carlo integration of a highly-oscillatory function, and show the advantage of using the density gradient in computing integrals of this type. In Section \ref{sec:general}, we extend all the concepts introduced in Section \ref{sec:1D} to higher-dimensional manifolds. Section \ref{sec:recursion} focuses on a recursive algorithm for the density gradient defined on a sequence of evolving manifolds under a differentiable map $\varphi$. Sections \ref{sec:1D}-\ref{sec:recursion} include examples of $x(\xi)$ defined by popular dynamical systems, as well as numerical results validating the derived expressions. Finally, Section \ref{sec:conclusion} concludes the paper.

%SECTION 2 ---------------------

\section{Computing $g$ on one-dimensional manifolds}\label{sec:1D}
In this section, we focus on the computation of the density gradient $g$ in the simplest topological setting. In particular, we consider one-dimensional manifolds, which can be described using a single parameter $\xi\in[0,1]$. That manifold is a curve $\mathcal{C}$ immersed in the Euclidean $\mathbb{R}^n$ space. We assume there exists a one-to-one map $x(\xi)\in\mathcal{C}\subset \mathbb{R}^n$, which is at least twice differentiable with respect to $\xi$, i.e., $x(\xi)\in C^2[0,1]$. In this case, the density gradient function is a scalar quantity defined as a directional derivative along $\mathcal{C}$ of logarithmic density, $g = \partial_s\log\rho$, where $\rho:\mathcal{C}\to[0,1]$ is a density function implied by $x(\xi)$. If we think of $\xi$ as a realization of the random variable uniformly distributed in $[0,1]$, then $x(\xi)$ is in fact the inverse cumulative distribution function (inverse CDF, a.k.a. the quantile function). Intuitively, $x(\xi)$ tells us that $100\xi\;\%$ of all points mapped from the uniformly distributed set are located on the curve segment between $x(0)$ and $x(\xi)$. On the other hand, the density function $\rho$ indicates the density of points mapped on $\mathcal{C}$ per unit curve length. Therefore, $\rho$ is counter-proportional to the magnitude of the first derivative of $x(\xi)$.

\begin{figure}
\centering
\includegraphics[width=1\textwidth]{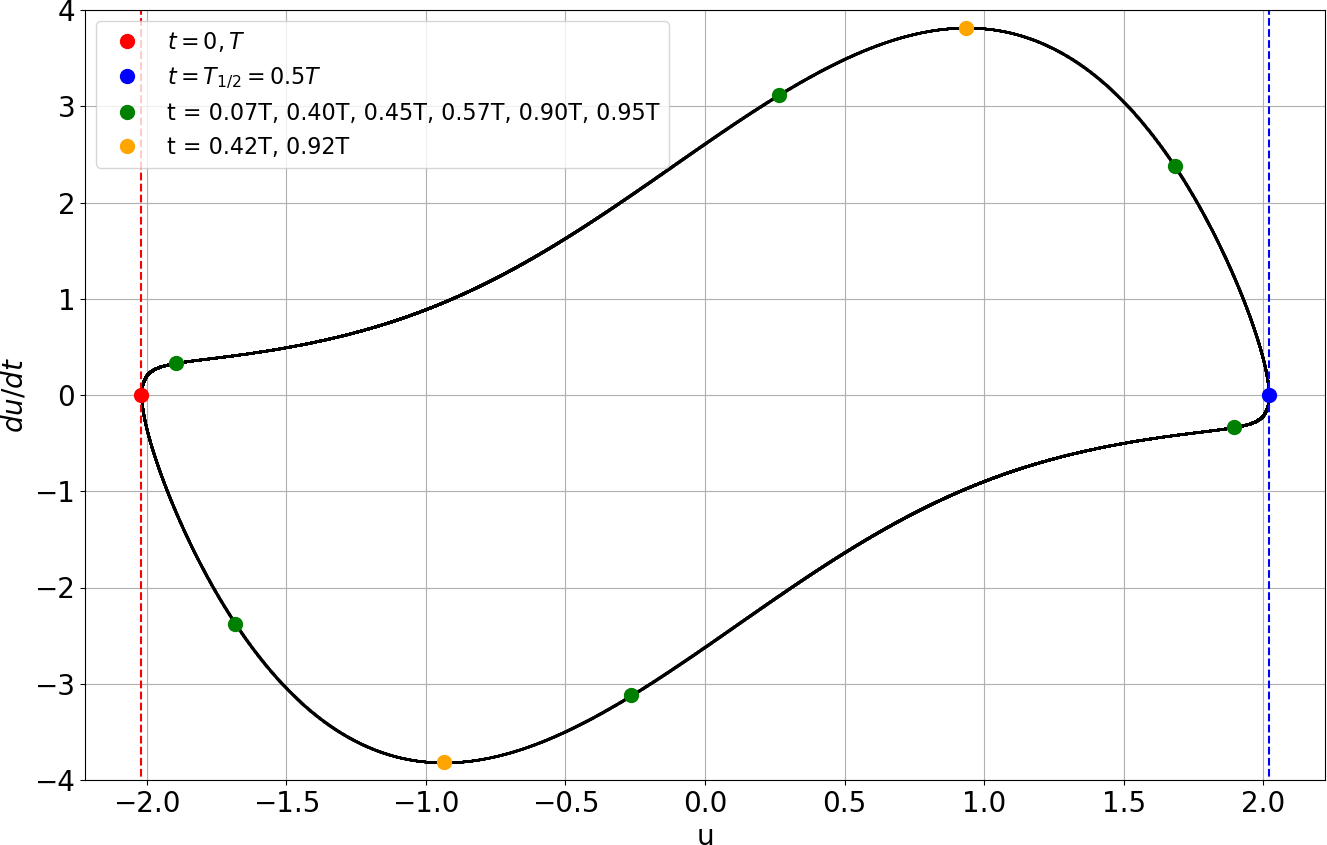}
\caption{Trajectory of the Van der Pol oscillator (Eq. \ref{eqn:van-der-pol}). The red dot represents the initial condition, as well as the solution after time $T$, while the blue dot indicates the solution after time $T_{1/2}$. The vertical dashed lines correspond to $u = -a$ and $u = a$ (boundaries of the range of $u$). At the green dots, the solution satisfies $du/dt+d^3u/dt^3 = 0$, while the zero acceleration state, $d^2u/dt^2=0$, is represented by orange dots.}
\label{fig:van-der-pol}
\end{figure}

In the following three subsections, we analytically derive the expression for $g$ in terms of the inverse CDF $x(\xi)$ for simple line manifolds, $n=1$ (Section \ref{sec:lines}), and general curves, $n\geq 1$ (Section \ref{sec:curves}), and demonstrate its importance in a numerical integration experiment (Section \ref{sec:integral}). We illustrate all relevant concepts using a certain $x(\xi)$ associated with the Van der Pol equation,
\begin{equation}
    \label{eqn:van-der-pol}
    \frac{d^2u}{dt^2} = 2(1-u^2)\frac{du}{dt} - u, \;\;u(0) = -a,\;\;\frac{du}{dt}(0) = 0,
\end{equation}
which describes the coordinates of a 2D non-conservative oscillator with non-linear dumping \cite{ginoux-vanderpol}. In our numerical examples, we choose $a = 2.0199$, in which case the solution $[u(t), du/dt(t)]^T$ approximately lies on the limit cycle with period $T = 2T_{1/2}\approx 7.638$ and $u(t)\in[-a,a]$ for all $t\geq 0$. Figure \ref{fig:van-der-pol} illustrates the limit cycle of Eq. \ref{eqn:van-der-pol}, which has been computed using the second-order Runge-Kutta (midpoint) method with time step $\Delta t = 0.0001$.

\subsection{Lines: $\mathcal{C}\subset\mathbb{R}$}\label{sec:lines}
We start from the simplest case, i.e., when $\mathcal{C}$ is a bounded line segment in $\mathbb{R}$, between $a$ and $b$. The corresponding inverse CDF $x(\xi)$ differentiably maps $[0,1]$ to $[a,b]$ and is related to the density function by the following expression,
\begin{equation}
\label{eqn:line1}
\xi(x) = \int_{a}^{x(\xi)} \rho(y)\;dy \;\;\;\forall x \in [a,b].    
\end{equation}
Since $\xi\in[0,1]$, $\rho(x)$ is in fact the probability density function (PDF) corresponding to the CDF $\xi(x)$, which satisfies $d\xi=\rho(x)\; dx$. Using the inverse function theorem, which asserts $f'(f^{-1}(c)) = 1/[(f^{-1})'(c)]$ for any differentiable one-to-one function $f$ at any $c$ such that $(f^{-1})'(c)\neq0$, we conclude that
\begin{equation}
\label{eqn:line2}
\frac{dx}{d\xi}(\xi)\;\rho(x(\xi)) = 1.   
\end{equation}
Eq. \ref{eqn:line2} indicates that at any point $x(\xi)$ on the manifold, the product of the PDF and derivative of the inverse CDF is constant. Thus, by differentiating Eq. \ref{eqn:line2} with respect to $\xi$ and reshuffling terms, we obtain a direct expression for $g$ at each point on the manifold,
\begin{equation}
\label{eqn:line3}
g(x(\xi)) = \partial_{x}\log\rho(x(\xi)) = \frac{\partial_x\rho(x(\xi))}{\rho(x(\xi))} = -\frac{\frac{d^2 x}{d\xi^2}(\xi)}{\Big(\frac{d x}{d\xi}(\xi)\Big)^2}.  
\end{equation}
To illustrate these functions and their relation, we will consider the solution to Eq. \ref{eqn:van-der-pol}, $u(t)$, for $t\in[0,T_{1/2}]$, where $T_{1/2} \approx 3.819$. Based on Figure \ref{fig:van-der-pol}, it is evident that $u(t)$ is a one-to-one smooth function and $du/dt \geq 0$ in that time interval. In fact, we can apply the linear transformation $t\to \xi$ to notice that
\begin{equation}
\label{eqn:line4}
 x(\xi) = u\left(\xi T_{1/2}\right)
\end{equation}
is a representation of the inverse CDF. Next, we compute the first and second derivative of Eq. \ref{eqn:line4} with respect to $\xi$ and plug them to Eq. \ref{eqn:line3} to obtain the following formula for $g$ along the trajectory,
\begin{equation}
    \label{eqn:line5}
    g(u(t)) = -\frac{\frac{d^2u}{dt}(t)}{\left(\frac{du}{dt}(t)\right)^2} \stackrel{\text{Eq. \ref{eqn:van-der-pol}}}{=} -\frac{2(1-u^2(t))\frac{du}{dt}(t) - u(t)}{\left(\frac{du}{dt}(t)\right)^2}.
\end{equation}
We observe that the density gradient is invariant to any linear change of variables, i.e., when $d\xi/dt$ is constant. Given a numerical solution to Eq. \ref{eqn:van-der-pol}, the density can be directly computed from $\rho(u(t)) = (T_{1/2}\; du/dt (u(t)))^{-1}$, which follows from Eq. \ref{eqn:line2}, whereas the density gradient function can be evaluated using Eq. \ref{eqn:line5}. 

\begin{figure}
\includegraphics[width=1.\textwidth]{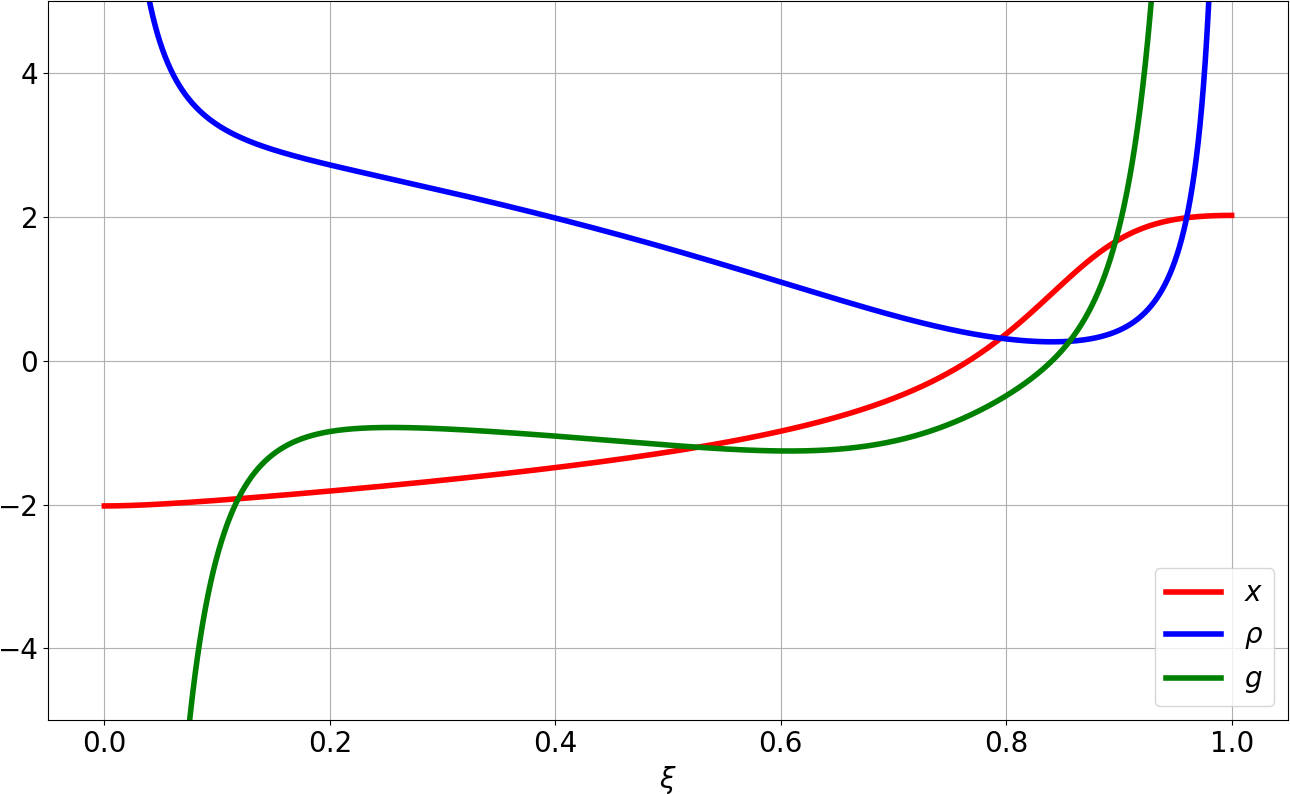}
\caption{The inverse CDF function $x(\xi)$ defined by the solution to the Van der Pol equation, such that $x(\xi(t)) = u(t)$ for all $t\in [0, T_{1/2}]$ (red), and the corresponding density (blue) and density gradient function (green). We used data presented in Figure \ref{fig:van-der-pol} to compute all the three functions.}
\label{fig:line}
\end{figure}

Figure \ref{fig:line} illustrates the inverse CDF $x(\xi)$, defined by Eq. \ref{eqn:line4}, as well as the corresponding density and density gradient. We clearly observe that both $\rho$ and $g$ are undefined at the endpoints, i.e., at $\xi = 0$ and $\xi = 1$, which is a consequence of zero slope of $x(\xi)$. Moreover, the larger the rate of change of $x$, the smaller the value of $\rho$, which confirms our previous intuitive explanation of the density function. We also notice that the density gradient is zero at the point corresponding to a local extremum of $\rho$ and the inflection point of $x(\xi)$.  

\subsection{Approximating integrals of a highly oscillatory function}\label{sec:integral}

We now demonstrate the use of the density gradient function in the numerical computation of a highly oscillatory function. Consider the following Lebesgue integral,
\begin{equation}
\label{eqn:line-integral}
    I = \int_{-a}^{a} \partial_xf(x)\;d\xi(x),
\end{equation}
where $\xi(x)$ denotes a Lebesgue measure defined by Eq. \ref{eqn:line1}, while $f$ is a function whose first derivative is integrable and bounded. Certainly, it is assumed the above integral converges. Indeed, a sufficient condition for the convergence of $I$ in this case is Lebesgue-integrability of the density gradient with respect to the density $\rho$ \cite{sliwiak-differentiability}, i.e., $g\in L^1(\rho)$.  However, the necessary and sufficient condition imposes extra requirements for the $f$ function itself, i.e., $\partial_x f\in L^1(\rho)$ or, equivalently, $\partial_x f\,\rho \in L^1[-a,a]$. In our experiment, the function $f$ has the following form,
\begin{equation}
    \label{eqn:line-function}
    f(x) = \left((x-a)(x+a)\sin(Kx^2)\right)^2,
\end{equation}
with some positive number $K$. We use Eq. \ref{eqn:line1} to rewrite the above integral, and then integrate it by parts. There exist a few scenarios when the resulting boundary term vanishes. One option is that the product $\partial_x f\,\rho$ is periodic and integrable on $[-a,a]$. Another possibility is when both $\partial_x f$ and $\rho$ are bounded and at least one of them vanishes at the domain boundaries. In any case, two new versions of $I$, alternative to the original form (in Eq. \ref{eqn:line-integral}), are available,
\begin{equation}
\label{eqn:line-integral2}
    \int_{-a}^{a} \partial_xf(x)\;\rho(x)\;dx = I =
    - \int_{-a}^{a} f(x)\;g(x)\;d\xi(x).
\end{equation}
To numerically approximate the integral $I$, we apply three distinct approaches. The integral in Eq. \ref{eqn:line-integral} and the RHS of Eq. \ref{eqn:line-integral2} can be estimated using a Monte Carlo method, which requires generating a random sequence $\{x^1,x^2,...,x^N\}$ distributed according to the measure $\xi$. If such a sequence is available, then the integral of any Lebesgue-integrable function $h(x)$ can be approximated as follows,
\begin{equation}
\label{eqn:line-integral3}
    \int_{-a}^{a} h(x)\;d\xi(x) \approx \frac{1}{N}\sum_{i=1}^N h(x^{i}),
\end{equation}
since $\xi\in[0,1]$. Finally, the integral on the LHS of Eq. \ref{eqn:line-integral2} is evaluated using a standard trapezoidal rule with a uniform $N$-element grid for $x$ between $-a$ and $a$. In the calculation, we allot
\begin{equation}
    \label{eqn:line-integral4}
x^{i}=u\left(\frac{i-1}{N-1}T_{1/2}\right).
\end{equation}
It can be numerically verified that for this particular choice of the sequence, $g\notin L^1(\rho)$, but $\partial_x f\in L^1(\rho)$. It means means that the integral $I$ converges despite the blow-up of $\rho$ and $g$ at the boundaries of $[-a,a]$. To assess the Lebesgue-integrability of these functions, we applied the procedure described in Section 4 of \cite{sliwiak-differentiability}. This algorithm approximates the slope of the distribution tail of any function in the logarithmic scale.  

\begin{figure}
\includegraphics[width=0.49\textwidth]{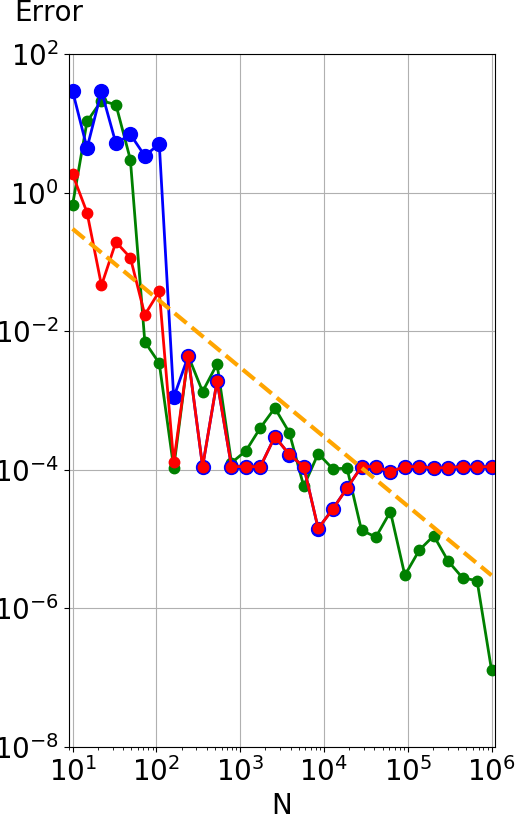}
\includegraphics[width=0.49\textwidth]{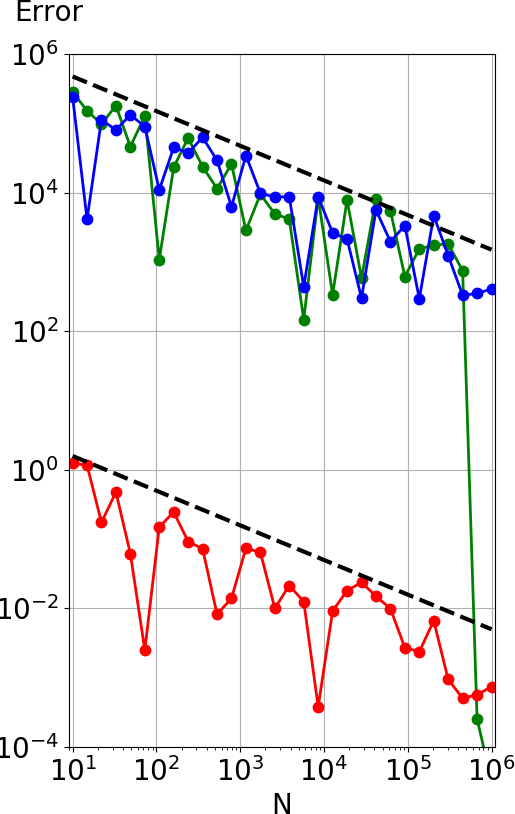}
\caption{Relative error of the approximation of $I$ for $K = 10$ (left) and $K = 100000$ (right) obtained using three methods: Monte Carlo integration applied to Eq. \ref{eqn:line-integral} (blue), Monte Carlo integration applied to the RHS of Eq. \ref{eqn:line-integral2} (red), and trapezoidal rule applied to the LHS of Eq. \ref{eqn:line-integral2} (green). Black and orange dashed lines are reference lines representing functions proportional to $N^{-1/2}$ and $N^{-1}$, respectively. In each of these plots, we computed the relative error with respect to the approximation of $I$ obtained using the trapezoidal with $N = 10^8$ samples.}
\label{fig:convergence}
\end{figure}

In order to compare the performance of these three integration methods, we proceed as follows. First, we generate the sequence $\{x^1,x^2,...,x^{N}\}$, $N = 10^5$ (time step is chosen such that $\Delta t = T_{1/2}/(N-1)$) and, using Eq. \ref{eqn:line2} and Eq. \ref{eqn:line5}, we directly evaluate $\rho$ and $g$ at all points from that sequence. Subsequently, both the density and density gradient functions are linearly interpolated everywhere between $-a$ and $a$. We use these interpolators to approximate the two functions at any point of a uniform grid (trapezoidal rule) or sequence defined by Eq. \ref{eqn:line-integral4} (Monte Carlo) for an arbitrary value of $N$. If $K$ is sufficiently small, then the approximation error of the trapezoidal rule is expected to be upperbounded by $\mathcal{O}(1/N)$, because the integrand, $\partial_x f\,\rho$, is Lebesgue-integrable \cite{cruz-traprule}. According to the Nyquist-Shannon sampling theorem, however, the discrete representation of the integrand may not be captured properly if $K$ is very large, in which case the trapezoidal rule's error decays as in a typical Monte Carlo method. Figure \ref{fig:convergence} shows the behavior of the relative error of the approximation of $I$ obtained using these three methods. The error is computed with respect to the reference solution obtained through the trapezoidal rule using $N=10^8$ points.       

We observe that for a moderately-oscillatory integrand ($K = 10$), the relative error of the trapezoidal rule (green curve) decays as $\mathcal{O}(1/N)$, which confirms the theoretical estimates. In this case, the performance of both the Monte Carlo approximations (blue and red curves) does not differ much from the trapezoidal rule's. The Monte Carlo approximation clearly converges to a solution slightly different than the reference solution, which is a consequence of the fact the latter was generated using the trapezoidal rule for a linearly interpolated function. This example indicates that there is no reason to perform integration by parts and compute $g$ to approximate integrals of  low- or  moderately-oscillatory functions. The right-hand side plot of Figure \ref{fig:convergence} corresponds to a different scenario, i.e., when $f$ is highly-oscillatory ($K=10^5$). The error now decays $\mathcal{O}(1/\sqrt{N})$ at $N\in[10^1,5\cdot10^5]$, regardless of the integration method. The trapezoidal rule requires almost $10^6$ samples to guarantee satisfactory accuracy. Note $\partial_x f$ has a magnitude proportional to $K$, and thus the variance of the sequence $\{\partial_x f(x^1),\partial_x f(x^2),...,\partial_x f(x^N)\}$ is of the order of $K^2$. Therefore, the Monte Carlo approach applied to Eq. \ref{eqn:line-integral} requires $\mathcal{O}(10^{10})$ samples to secure error of the order of 1. A similar error can be achieved if we perform integration by parts and compute $g$ and generate only $\mathcal{O}(1)$ samples, since the variance is reduced $10^{10}$ times. 

In conclusion, the computational cost of the Monte Carlo method can be dramatically reduced using the generalized fundamental calculus theorem. In case of the $f$ function, the regularization of the integral in Eq. \ref{eqn:line-integral} may decrease the cost even $K^2$ times. This result is significant specifically in the context of strongly fluctuating functions.

\subsection{General curves: $\mathcal{C}\subset \mathbb{R}^{n}$}\label{sec:curves}

We extend the concepts introduced in Section \ref{sec:lines} to the case in which $x(\xi)$ differentiably maps $[0,1]$ to $\mathcal{C}\subset \mathbb{R}^{n}$, where $n$ is some positive integer. Geometrically, $x(\xi)$ represents a curve embedded in the $n$-dimensional Euclidean space. The measure $\xi(x)$ can now be expressed as an integral of the density, $\rho: \mathcal{C}\to [0,1]$, along $\mathcal{C}$ with respect to the arc length $s$,
\begin{equation}
\label{eqn:curve1}
    \xi(x) = \int_{\mathcal{C}[x(0),\,x(\xi)]} \rho(x)\;ds,
\end{equation}
where $\mathcal{C}[x(0),x(\xi)]$ denotes a segment of $\mathcal{C}$ between the points indicated in the square bracket. Due to the parameterization $x(\xi)$, the length of the curve $\mathcal{C}$ equals $\int_{\mathcal{C}}ds$, while the arc length differential $ds$ is related to $d\xi$ by $ds = \|dx/d\xi\|\;d\xi$. Using this relation and Eq. \ref{eqn:curve1}, we obtain the following identity,
\begin{equation}
\label{eqn:curve2}
    \rho(x(\xi))\;\left\|\frac{dx}{d\xi}(\xi)\right\| = 1.
\end{equation}
We now differentiate Eq. \ref{eqn:curve2} with respect to $\xi$, apply the chain rule and reshuffle terms,
\begin{equation}
\label{eqn:curve3}
    g(x(\xi))= \partial_s\log(\rho(x(\xi))) =\frac{\partial_{s}\rho}{\rho}(x(\xi)) = -\frac{\frac{d x}{d\xi }(\xi)\cdot \frac{d^2 x}{d\xi^2}(\xi)}{\| \frac{d x}{d\xi }(\xi)\|^3},
\end{equation}
where $\partial_s$ denotes the directional derivative along the curve $\mathcal{C}$ in the direction of increasing $\xi$. Note the expression for $g$ in Eq. \ref{eqn:curve3} reduces to Eq. \ref{eqn:line3} if $x(\xi)$ represents a line manifold, i.e., $\mathcal{C}\subset\mathbb{R}^1$.

\begin{figure}
\includegraphics[width=1.05\textwidth]{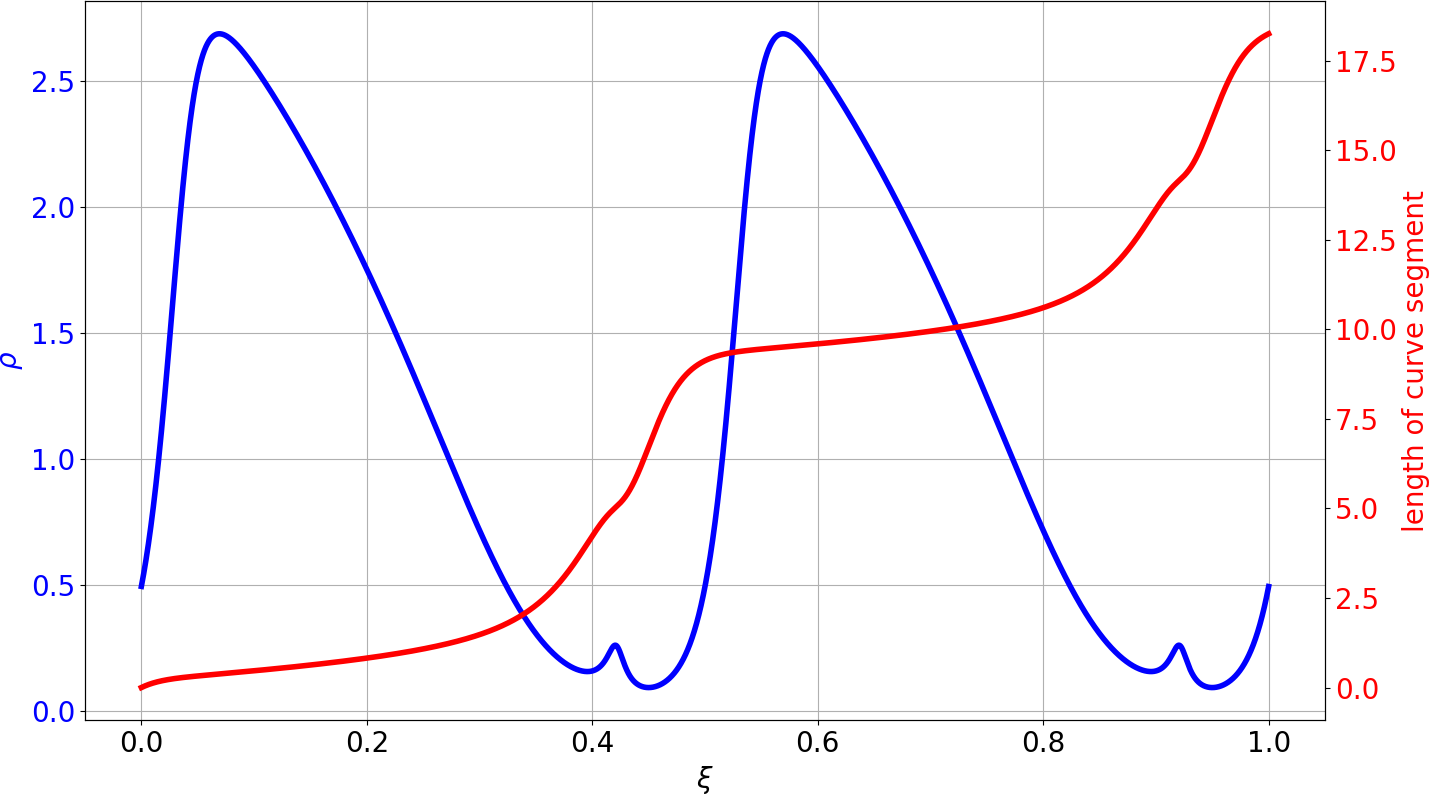}
\caption{The density $\rho$ (blue) and length of the curve segment $\mathcal{C}[x(0),x(\xi)]$ (red) associated with the map $x(\xi)$ defined by Eq. \ref{eqn:curve4}. The former is computed using the analytical expression in Eq. \ref{eqn:curve2}, while the latter is approximated by summing the length of consecutive linear segments connecting the points in the sequence $\{x(0), x(\Delta t/T), x(2\Delta t/T),...,x(\xi)\}$, obtained in the numerical integration of Eq. \ref{eqn:van-der-pol}.}
\label{fig:curve-density}
\end{figure}

\begin{figure}
\includegraphics[width=1.\textwidth]{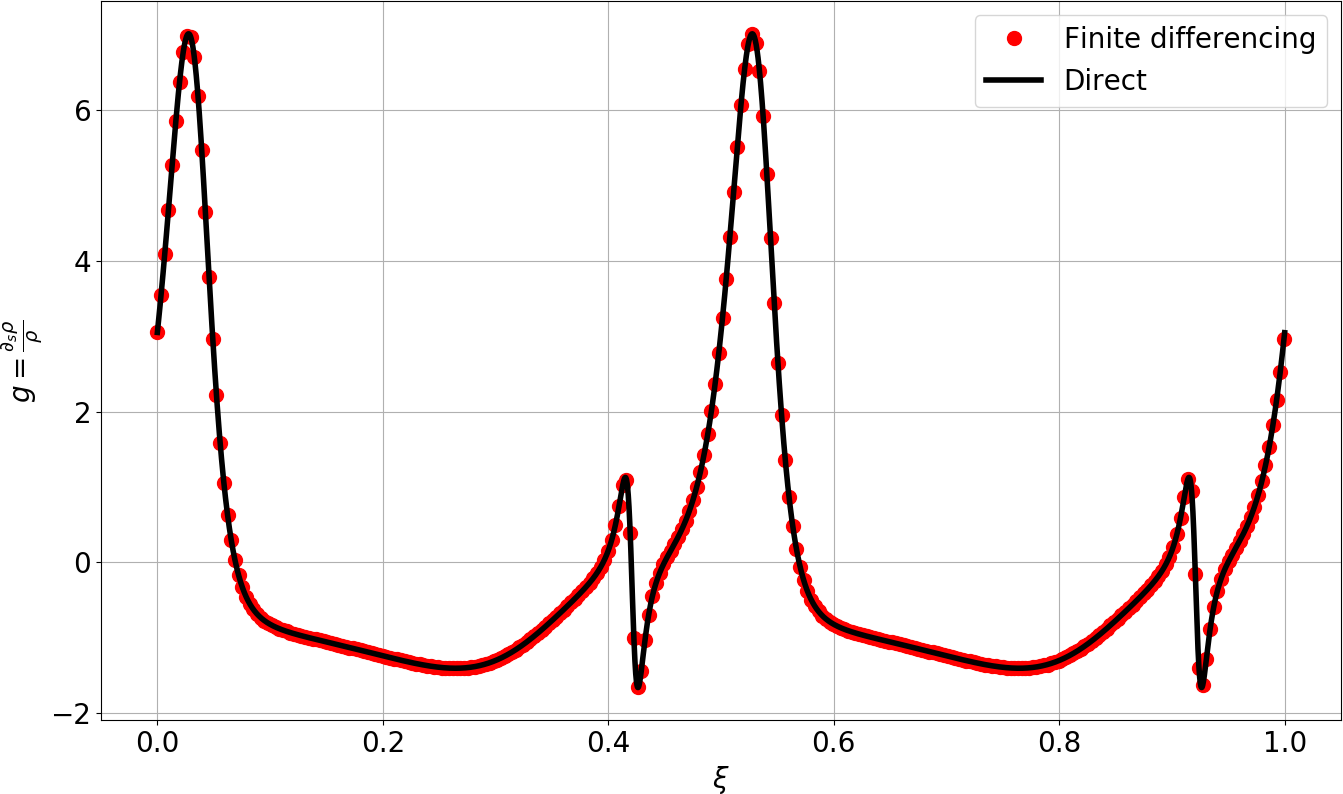}
\caption{The density gradient function $g$ computed directly (using Eq. \ref{eqn:curve3}) and through a finite difference method. In the latter approach, we note $\partial_s \rho = \partial_{\xi}\rho / \partial_{\xi} s$. Both the numerator and denominator is approximated using the central finite difference scheme on a uniform grid using data presented in Figure \ref{fig:curve-density}.}
\label{fig:curve-gradient}
\end{figure}

As an example, we re-consider the Van der Pol oscillator (Eq. \ref{eqn:van-der-pol}). This time, however, $x(\xi)$ represents a curve embedded in $\mathbb{R}^2$. In particular, $x(\xi)$ describes a two-dimensional loop such that
\begin{equation}
    \label{eqn:curve4}
    x(\xi) = \begin{bmatrix}u(\xi T)\\\frac{du}{dt}(\xi T)\end{bmatrix}
\end{equation}
(see Figure \ref{fig:van-der-pol} for an illustration of the loop). If a numerical solution to Eq. \ref{eqn:van-der-pol} is available, one can combine Eq. \ref{eqn:curve2} with Eq. \ref{eqn:curve4} to directly evaluate the density function. Similarly, by plugging Eq. \ref{eqn:curve4} to Eq. \ref{eqn:curve3}, it is possible to compute the density gradient, analogously to the procedure described in Section \ref{sec:lines}. Consequently, on the RHS of Eq. \ref{eqn:curve3}, $dx/d\xi$ can be replaced with $du/dt$, and $d^2x/d\xi^2$ with $d^2u/dt^2$. We can do so because the density gradient is invariant to any linear transformation of variables. 
Figure \ref{fig:curve-density} illustrates the density function $\rho$, as well as the length of the curve segment $ \mathcal{C}[x(0),x(\xi)]$, versus the parameter $\xi$. We observe $\rho$ is large if the slope of the length function is small, and vice versa, which is analogous to the $x\--\rho$ relation in Figure \ref{fig:line}. In this particular case, $\rho(\xi)$ is clearly a periodic function with period 0.5. This property is manifested in Figure \ref{fig:van-der-pol}. Indeed, one can notice the relation between $du/dt$ and $u$ at $t\in[0,T_{1/2}]$ is the same as $-du/dt$ and $-u$ at $t\in[T_{1/2},T]$, where $T_{1/2}$ corresponds to $\xi = 0.5$. 
Figure \ref{fig:curve-gradient} shows the density gradient $g$ computed using two distict ways: through a direct evaluation via Eq. \ref{eqn:curve3} and a finite difference scheme (see the caption for more details). The two approaches provide visibly identical solutions, which confirms the correctness of Eq. \ref{eqn:curve3}. Clearly, the density gradient inherits the periodic behavior of $\rho$. We notice that the density gradient vanishes if the numerator of Eq. \ref{eqn:curve3} is zero, which can happen if $d^2u/dt^2 = 0$ (at the two orange dots in Figure \ref{fig:van-der-pol}) and/or $du/dt + d^3u/dt^3 = 0$ (at the six green dots in Figure \ref{fig:van-der-pol}). These two cases coincide with the local extrema of the density function (i.e., $d\rho/d\xi = 0$ if at least one of these equations is satisfied). However, zero density gradient does not imply the inflection point ($d^2u/dt^2 = 0$), in contrast to the line manifold case (see Section \ref{sec:lines}).

%SECTION 3----------------------
\section{Computing $g$ on general smooth manifolds}\label{sec:general}
The purpose of this section is to generalize the concept of the density gradient and derive a formula for $g$ defined on higher-dimensional manifolds. Here, we consider a smooth invertible map $x(\xi):U\to M$, where $U\subset \mathbb{R}^m$, $M\subset \mathbb{R}^n$, $m\leq n$, $x = [x_1,...,x_n]^T$ and $\xi = [\xi_1,...,\xi_m]^T$. $U$ is an $m$-orthotope (hyperrectangle), which is defined as the Cartesian product of $m$ 1D line manifolds (i.e., intervals of the real line). We no longer assume that these elementary sets only involve numbers between 0 and 1. $M$ is an oriented differentiable manifold, whose shape is defined by the chart $x$. For example, if $m=2$ and $n=3$, then $M$ represents a smooth surface. The density gradient $g$ is now defined as a directional gradient of the logarithm of the density function $\rho:M\to[0,1]$ implied by the chart $x(\xi)$. In particular, $g = \nabla_s\log\rho$, $\nabla_s := [\partial_{s_1}, \partial_{s_2}, ..., \partial_{s_m}]^T$, where $\partial_{s_i}$, $i=1,...,m$, denote directional derivatives along the corresponding {\em isoparametric curves}. The $i$-th component of $g$ is the rate of change of $\log\rho$ along the curve whose preimage involves vectors $\xi\in U$ with constant all coordinates except $\xi_i$. If $\log\rho$ is differentiable with respect to all the coordinates of $x$ and $\nabla_x :=[\partial_{x_1}, \partial_{x_2}, ..., \partial_{x_n}]^T$, then $\partial_{s_{i}} \log\rho = \nabla_x\log\rho\cdot s_i$, where $s_{i}$ denotes the unit vector that is tangent to the corresponding isoparametric curve and points in the direction of increasing $\xi_i$. In Section \ref{sec:general-formula}, we derive a generic formula for $g$, while Section \ref{sec:general-example} provides a specific example of a two-dimensional smooth manifold embedded in $\mathbb{R}^3$ (with $m=2$ and $n=3$).      

\subsection{Derivation of the general formula}\label{sec:general-formula}
Recall $x(\xi):U\to M$ is an invertible and differentiable map, where $U\subset R^m$, $M\subset R^n$, and $m\leq n$, while $\rho(x):M\to [0,1]$ is the density function implied by that chart. Let $\omega(x)$ be the natural volume form defined on $M$. Therefore, the Lebesgue measure $m$ of any subset $V\subset U$, mapped by $x$ to $N\subset{M}$, equals
\begin{equation}
m(V) = \int_{N}\rho(x)\;d\omega(x),
\end{equation}
which implies that the volume element $dm$ defined on $U$ can be expressed in terms of $\rho$ and the volume element defined on $M$, at every point $x(\xi)$,
\begin{equation}
\label{eqn:general-measure}
dm = d\xi_1\wedge d\xi_2\wedge ... \wedge d\xi_m = \rho(x)\;d\omega(x) = \rho(x)\;dx_1\wedge dx_2\wedge ... \wedge dx_n.
\end{equation}
The wedge symbol ($\wedge$) denotes the exterior product, while $d\xi_{i}$, $i=1,...,m$ and $dx_{i}$, $i=1,...,n$ represent covectors (1-forms) associated with the corresponding coordinate directions. Intuitively, these 1-forms measure small displacements in the direction of one coordinate. The volume element on $M$, $d\omega$, can be expressed in terms of $\xi$,
\begin{equation}
\label{eqn:general-substitution}
d\omega(x(\xi)) = \sqrt{\det C(x(\xi))}\;d\xi_1\wedge d\xi_2\wedge ... \wedge d\xi_m,
\end{equation}
where $C$ represents the $m\times m$ metric tensor of the coordinate transformation $\xi \to x$, defined as
\begin{equation}
\label{eqn:general-metric}
C(x(\xi)) = [\nabla_{\xi} x(\xi)]^T\;\nabla_{\xi} x(\xi),
\end{equation}
or, componentwise,
\begin{equation}
C_{ij}(x(\xi)) = \partial_{\xi_i} x(\xi) \cdot \partial_{\xi_j} x(\xi).
\end{equation}
The vector gradient $\nabla_{\xi} x(\xi)$ is represented by an $n\times m$ matrix, in which the $j$-th column contains the derivative of $x$ with respect to $\xi_j$, i.e., $[\nabla_{\xi}x(\xi)]_{ij} = \partial_{\xi_{j}}x_{i}(\xi)$.  
Combining Eq. \ref{eqn:general-measure} and \ref{eqn:general-substitution}, we conclude that the relation between the density function $\rho$ and metric tensor $C$, at any point $x(\xi)\in{M}$, can be written in the following way,
\begin{equation}
    \label{eqn:general1}
    \rho(x(\xi))\;\sqrt{\det C(x(\xi))} = 1,
\end{equation}
which is a generalization of Eq. \ref{eqn:curve2}. 
Let us now QR-factorize the vector gradient $\nabla_{\xi} x(\xi)$,
\begin{equation}
    \label{eqn:general-qr}
    \nabla_{\xi} x(\xi) = Q(x(\xi))\;R(x(\xi)),
\end{equation}
where $Q$ is an $n\times m$ matrix, whose columns form an orthonormal basis for the column space of $\nabla_{\xi} x(\xi)$, while $R$ is an $m\times m$ upper-triangular matrix. Note $Q^TQ = I$ everywhere on $M$. Using this property, we immediately notice that $C = R^T R$ and, therefore, Eq. \ref{eqn:general1} reduces to
\begin{equation}
    \label{eqn:general2}
    \rho(x(\xi))\;|\det R(x(\xi))| = 1.
\end{equation}
For any invertible matrix $A(s)$, which depends on a scalar $s$, the following indentity is true,
\begin{equation}
    \label{eqn:general-indentity}
    \frac{\partial \det A(s)}{\partial s} = \det A\; \mathrm{tr}\left(A^{-1}(s)\;\frac{\partial A(s)}{\partial s}\right). 
\end{equation}
Differentiating Eq. \ref{eqn:general2} with respect to $\xi_{i}$, applying chain rule and Eq. \ref{eqn:general-indentity}, we obtain the following expression for the $i$-th component of the density gradient, 
\begin{equation}
    \label{eqn:general-g1}
    g_i(x(\xi)) = \frac{\partial_{s_i}\rho(x(\xi))}{\rho(x(\xi))} = - \frac{\partial_{s_i} \det R(x(\xi))}{\det R(x(\xi))} = - \frac{\partial_{\xi_i} \det R(x(\xi))}{\det R(x(\xi))\|\partial_{\xi_i}x(\xi)\|}.
\end{equation}
Eq. \ref{eqn:general-g1} is computationally inconvenient, as it involves evaluating the determinant of $R$ and its directional derivative. Our goal is to rewrite the RHS of that equation such that only first and second parametric derivatives of $x(\xi)$, as well as $Q$ and $R$ factors, are involved.

Since $R$ is an upper-triangular matrix, we notice that
\begin{equation}
    \label{eqn:general-R}
    \frac{\partial \det R}{\det R} = \frac{\partial \left(\prod_{k=1}^m R_{kk}\right)}{\prod_{k=1}^m R_{kk}} = \sum_{k=1}^{m}\frac{(\partial R)_{kk}}{R_{kk}} = \mathrm{tr}(\partial R\;R^{-1}).
\end{equation}
Now, differentiating Eq. \ref{eqn:general-qr} with respect to $\xi_i$, and then left- and right-multiplying the resulting expression by $Q^T$ and $R^{-1}$, respectively, we obtain
\begin{equation}
    \label{eqn:general-aux1}
    Q^T(x(\xi)) \;\partial_{\xi_i}\nabla_{\xi}x(\xi)\;R^{-1}(x(\xi)) = Q^T(x(\xi))\;\partial_{\xi_i}Q(x(\xi)) + \partial_{\xi_i}R(x(\xi))\;R^{-1}(x(\xi)). 
\end{equation}
Note that since $Q^TQ = I$, then $Q^T\;\partial_{\xi_i} Q$ is anti-symmetric, which means its trace vanishes. Therefore, the following equality
\begin{equation}
    \label{eqn:general-treq}
    \mathrm{tr}\left(Q^T(x(\xi)) \;\partial_{\xi_i}\nabla_{\xi}x(\xi)\;R^{-1}(x(\xi))\right) = \mathrm{tr}\left(\partial_{\xi_i}R(x(\xi))\;R^{-1}(x(\xi))\right)
\end{equation}
holds everywhere on $M$. Finally, by combining Eq. \ref{eqn:general-g1}, \ref{eqn:general-R} and \ref{eqn:general-treq}, we obtain the general formula for $g_i$,
\begin{equation}
    \label{eqn:general-g2}
    g_i(x(\xi)) = \partial_{s_{i}}\log\rho(x(\xi)) = -\frac{\mathrm{tr}\left(Q^T(x(\xi)) \;\partial_{\xi_i}\nabla_{\xi}x(\xi)\;R^{-1}(x(\xi))\right)}{\|\partial_{\xi_i}x(\xi)\|},
\end{equation}
which holds everywhere on $M$ for $i = 1,...,m$. Using Einstein's summation convention, Eq. \ref{eqn:general-g2} can be rewritten to
\begin{equation}
\label{eqn:general-g2-einstein}
g_i(x(\xi)) = -\frac{q_j(x(\xi))\cdot\partial_{
\xi_i} \partial_{\xi_k} x(\xi)\;R_{kj}^{-1}(x(\xi))}{\|\partial_{\xi_i}x(\xi)\|},
\end{equation}
where $q_{j}(x(\xi))$ denotes the $j$-th column of $Q(x(\xi))$. Thus, to directly compute the density gradient at any point on a manifold, all first and second derivatives of the chart $x(\xi)$ must be found. In addition, QR factorization of the vector gradient $\nabla_\xi x$ and inversion of the $R$ matrix must be performed. In practice, inverting the triangular matrix $R$ means solving a linear system using the backward substitution method, which requires $\mathcal{O}(m^2)$ operations. Note Eq. \ref{eqn:general-g2-einstein} reduces to Eq. \ref{eqn:curve3} if $m=1$. In the following section, we present an example illustrating some of these quantities. Although Eq. \ref{eqn:general-g2-einstein} is a formula for the  derivative in the direction of a isoparametric curve, we can compute derivatives of $\log\rho$ in an arbitrary direction using the distributive law of the dot product. 

\subsection{Example: a surface manifold}\label{sec:general-example}
As an example of a surface manifold (with $m=2$ and $n=3$), let us consider $x(\xi) = u(\xi) =  [u_1(\xi),u_2(\xi),u_3(\xi)]^T$, where $\xi= [c,t]^T$, $-5 \leq c \leq 5$, $0 \leq t \leq 0.4$, $\left.u(\xi)\right|_{t=0} = [c,c,28]^T$, and $\partial_t u(\xi) = f(u(\xi))$, where $f$ is defined as follows,
\begin{equation}
    \label{eqn:general-lorenz}
    \begin{split}
     & \partial_t u_1(\xi) = 10\;(u_2(\xi)-u_1(\xi)),\\& \partial_t u_2(\xi) = u_1(\xi)\,(28 - u_3(\xi)) - u_2(\xi), \\ & \partial_t u_3(\xi) = u_1(\xi)\,u_2(\xi) - \frac{8}{3}u_3(\xi).
    \end{split}
\end{equation}
System \ref{eqn:general-lorenz} represents the Lorenz '63 oscillator, which is a mathematical model used for atmospheric convection \cite{lorenz-climate}. This system is known to exhibit chaotic behavior. However, we are interested in the solution in a short time interval, such that the trajectories do not intersect and the resulting surface is orientable. In particular, we compute $x(\xi)$ by numerically integrating System \ref{eqn:general-lorenz} in time for different values of $c\in[-5,5]$, using the second-order Runge-Kutta scheme with $\Delta t = 0.002$. There are two reasons we have chosen this particular $x(\xi)$. First, it serves as a perfect example of a problem, in which the smooth one-to-one solution, $x(\xi)$, cannot be found analytically. Thus, the computation of $g$ should be performed numerically using closed-form relations derived in Section \ref{sec:general-formula}. Second, the surface described by the chart $x(\xi)$ can be obtained as a evolution of 1D manifolds. This observation is utilized in Section \ref{sec:recursion}, where we derive expressions for evolving manifolds. To evaluate $\rho$ and $g$, we directly use Eq. \ref{eqn:general2} and Eq. \ref{eqn:general-g2}, respectively. To find these quantities, the vector gradient $\nabla_{\xi}x(\xi)=[\partial_c x(\xi), \partial_t x(\xi)]$, as well as the following second derivatives: $\partial_c^2 x(\xi), \partial_t^2 x(\xi), \partial_c\partial_t x(\xi)$, must be found at every point on the manifold. The time derivative, $\partial_t x(\xi) = f(x(\xi))$, is obtained automatically as we integrate System \ref{eqn:general-lorenz} in time. The second derivative of $x$ with respect to $t$ is obtained using the chain rule, $\partial_t^2 x(\xi) = \partial_t f(x(\xi)) = Df(x(\xi))\,f(x(\xi))$, where $Df$ denotes the Jacobian of System \ref{eqn:general-lorenz}. Thus, from the computational point of view, we need to solve a tangent equation to find  $\partial_t x(\xi)$ at every point of the trajectory defined by System \ref{eqn:general-lorenz}. Using this approach, one can analogously find derivatives with respect to $c$. Let $v(\xi) = \partial_c x(\xi)$ and $w(\xi) = \partial^2_c x(\xi)$. Using the chain rule, we conclude that $\partial_t v(\xi) = Df(x(\xi))\,v(\xi)$, $\left.v(\xi)\right|_{t=0} = [1,1,0]^T$ and, by differentiating again, $\partial_t w(\xi) = D^2f(x(\xi))(w(\xi), w(\xi))  + Df(x(\xi))\,w(\xi)$, $\left.w(\xi)\right|_{t=0} = [0,0,0]^T$, where $D^2 f$ denotes the Hessian of $f$. Using Einstein's summation convention, the $i$-th component of the blinear form $D^2f(x(\xi))(w(\xi), w(\xi))$ can be written as  $\partial_{x_k}\partial_{x_{l}}f_i\,w_k\,w_l$. Finally, the mixed derivative $\partial_c\partial_t x(\xi) = \partial_t v(\xi)$ is a byproduct of the numerical integration of the tangent equation for $v$. We solve all of these tangent equations using the same time integrator as the one mentioned above. Since $m=2$, the $2\times 2$ $R$ matrix is inverted analytically at every point on the trajectory.

\begin{figure}
\includegraphics[width=1.0\textwidth]{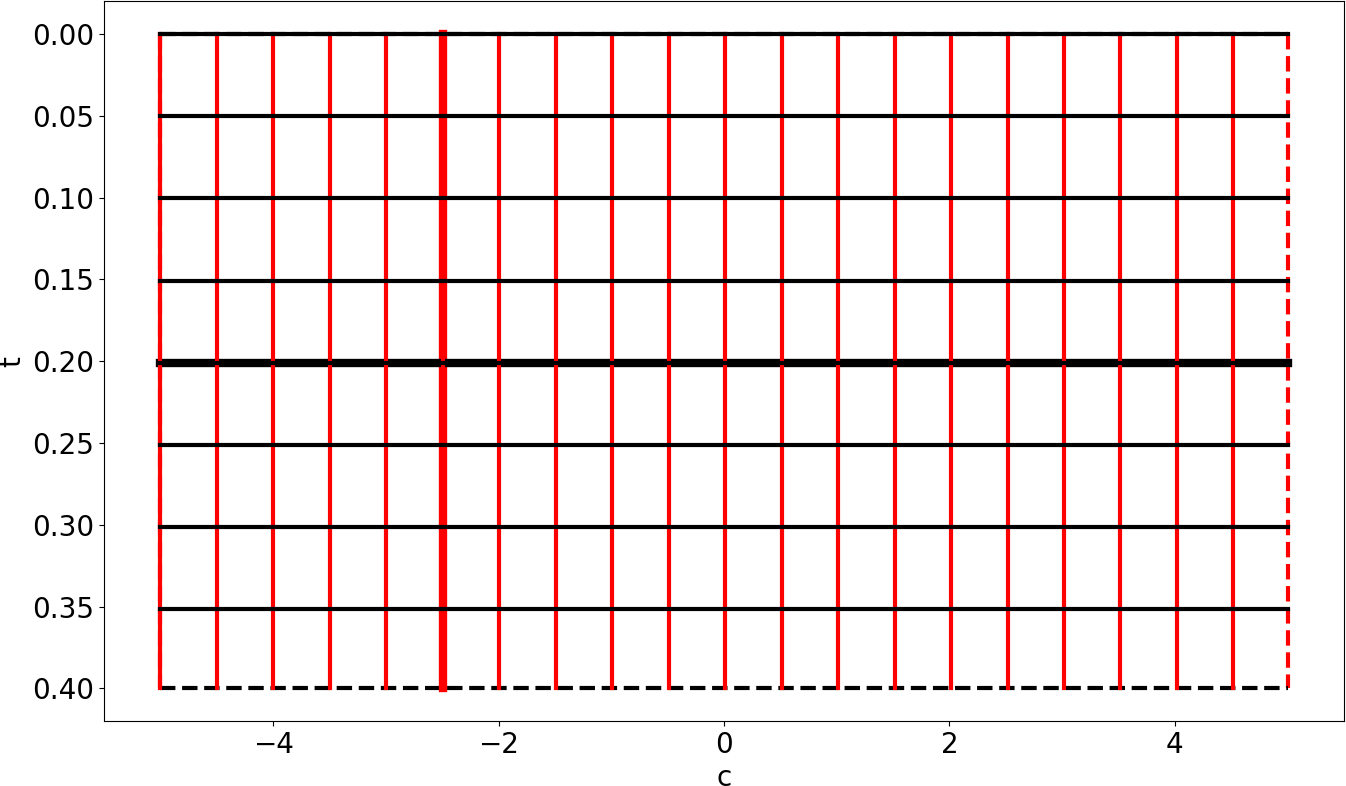}
\includegraphics[width=1.0\textwidth]{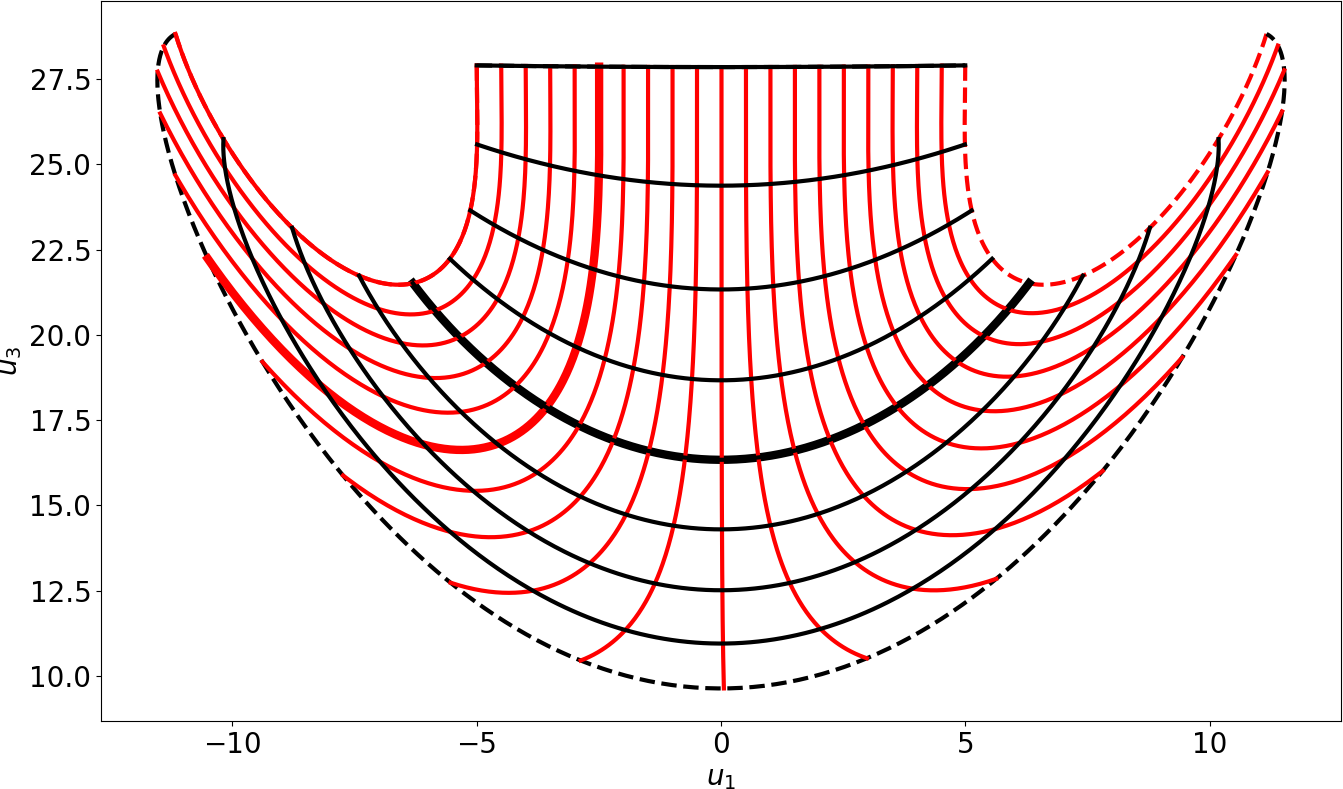}
\caption{Upper plot: a structured mesh representing the domain $U=\{(c, t)\,|\,c\in[-5,5],\, t\in[0,0.4]\}$. The black lines correspond to fixed values of $t$, while the red lines illustrate $\xi$ with a fixed value of $c$. The red and black dashed lines represent $c=5$ and $t=0.4$, while the red and black bold lines refer to $c = -2.5$ and $t=0.2$, respectively. Lower plot: $u_1\--u_3$ projection of the image of the structured mesh obtained through the mapping $x(\xi)$.}
\label{fig:general-mesh}
\end{figure}

\begin{figure}
\includegraphics[width=1.0\textwidth]{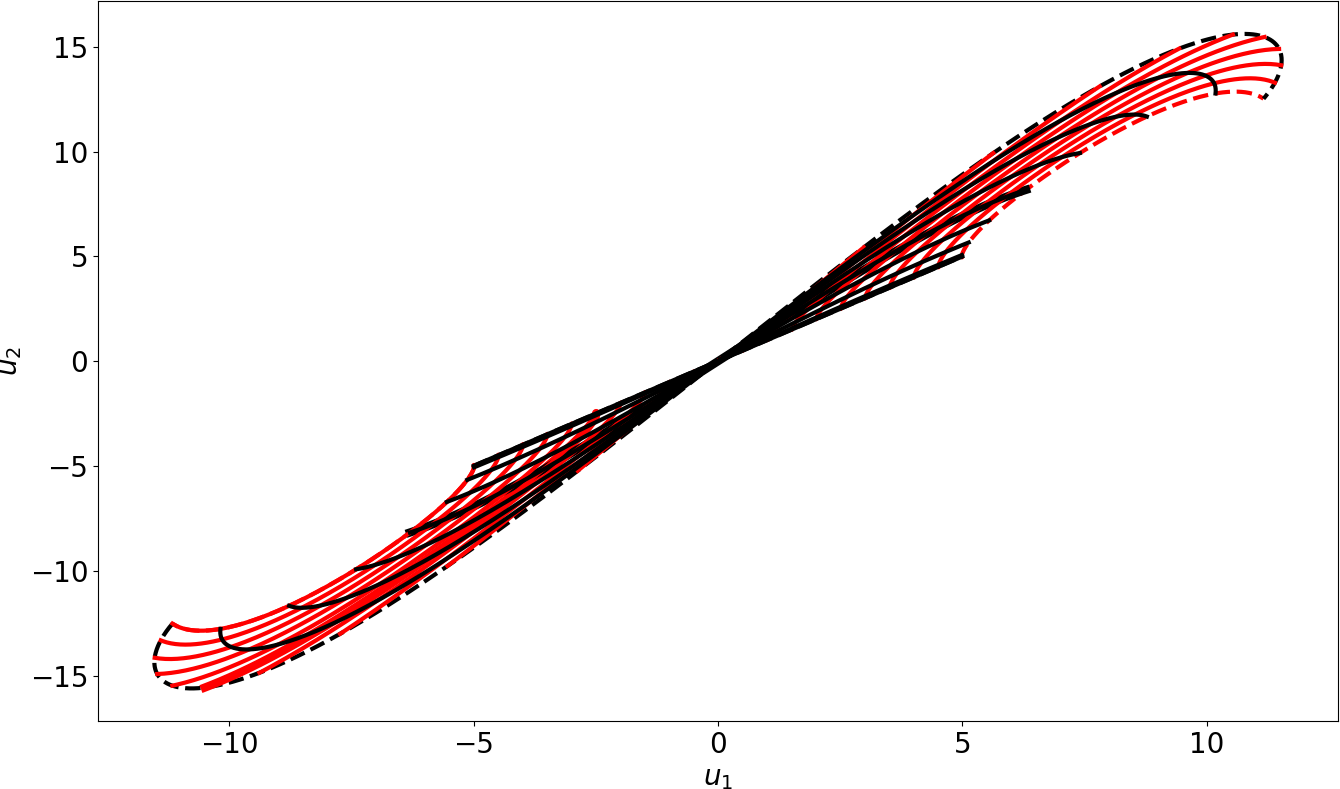}
\includegraphics[width=1.0\textwidth]{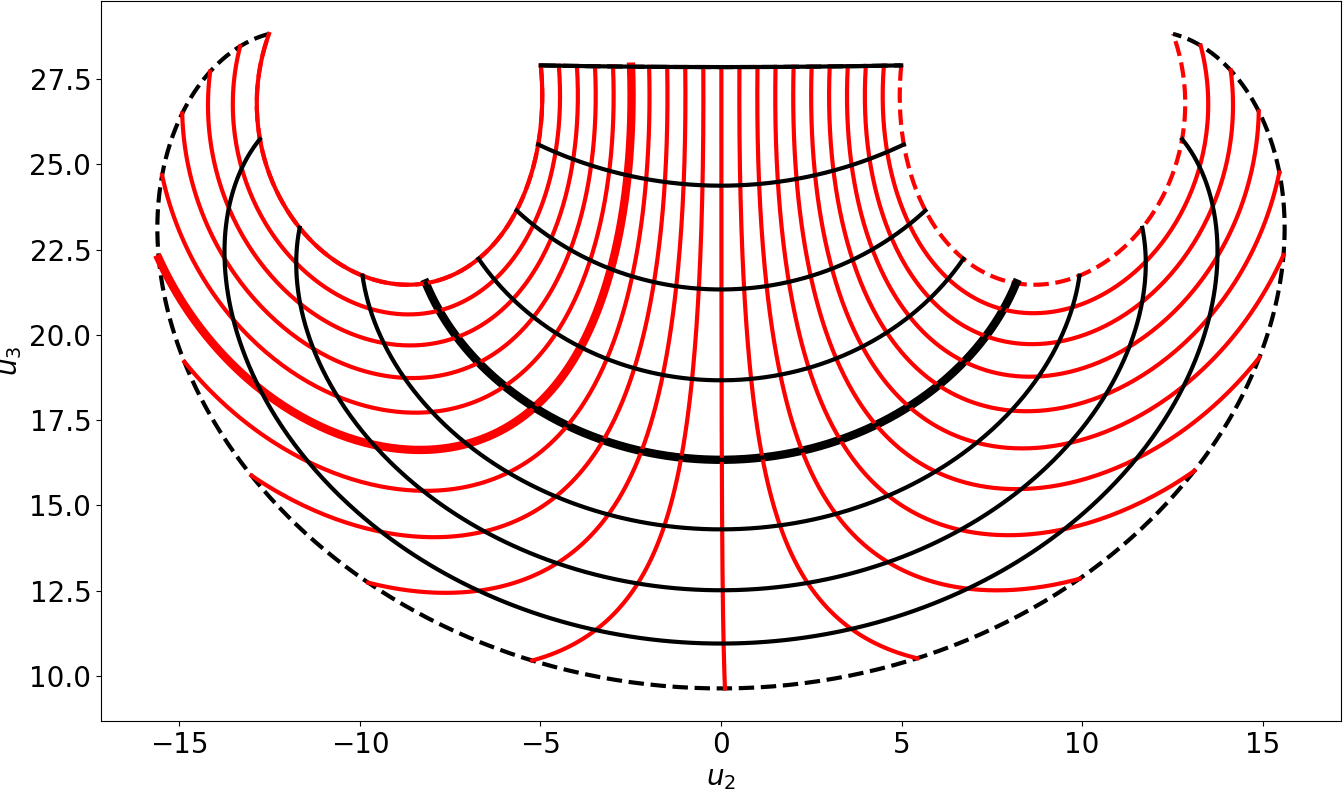}
\caption{Extension of Figure \ref{fig:general-mesh}. $u_1\--u_2$ (upper plot) and $u_2\--u_3$ (lower plot) projection of the image of the structured mesh obtained through the mapping $x(\xi)$.}
\label{fig:general-mesh-ext}
\end{figure}

\begin{figure}
\includegraphics[width=1.\textwidth]{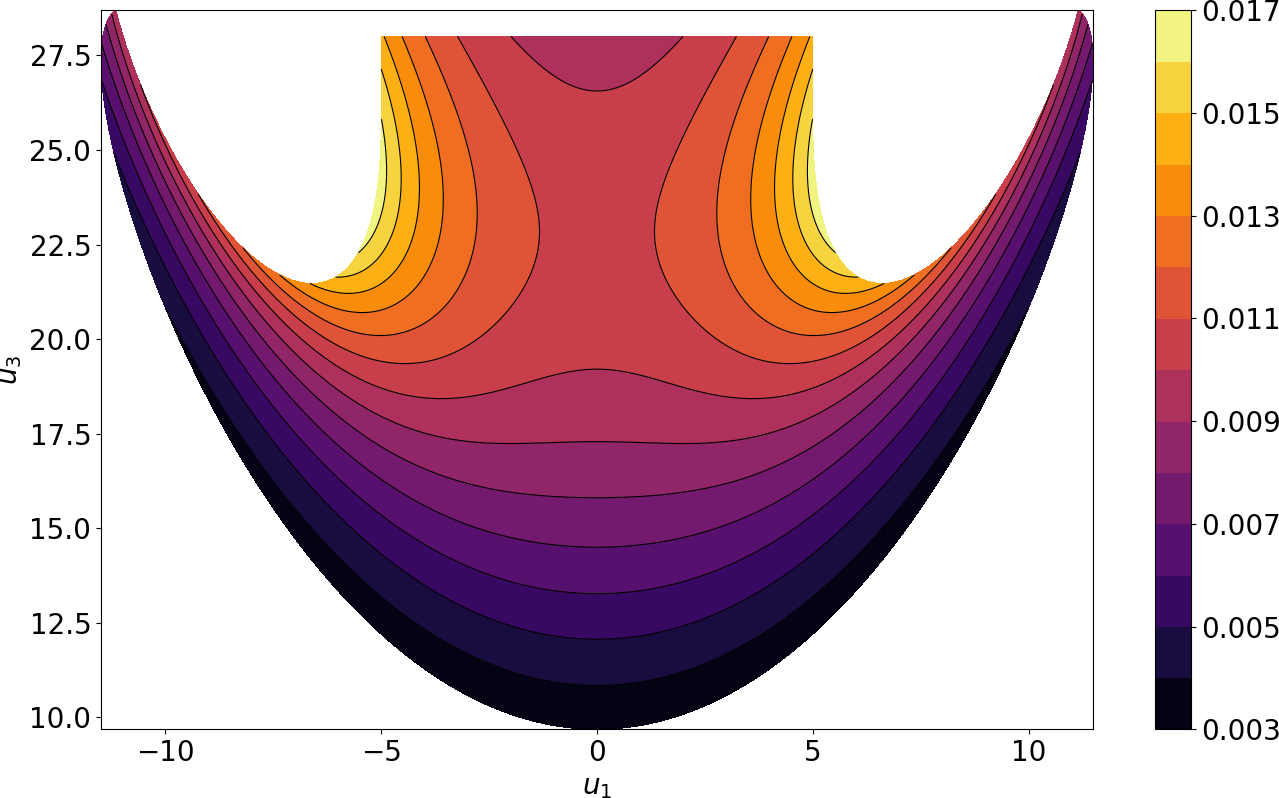}

\caption{$u_1\--u_3$ projection of the density function $\rho$.}
\label{fig:general-density}
\end{figure}

In this case, the $U$ space, which is the domain (preimage) of $x$, is in fact a Cartesian product of $[-5,5]$ and $[0,0.4]$. The upper plot in Figure \ref{fig:general-mesh} graphically represents $U$, while the lower plot illustrates the $u_1\--u_3$ projection of $M$, obtained through the mapping $x(\xi)$. For completeness, in Figure \ref{fig:general-mesh-ext}, we also include the $u_1\--u_2$ and $u_2\--u_3$ projection of the deformed mesh. It is clear that the deformation is symmetric with respect to $c = 0$. We also observe that fibers (isoparametric lines) corresponding to larger values of $t$ are subject to greater stretching than those at smaller $t$. These features are reflected by the distribution of the density function $\rho$, plotted in Figure \ref{fig:general-density}. The smaller the area of each distorted quadrilateral of the mesh, the larger the value of the density function. Indeed, the smallest values of the density distribution are located around $t=0.4$. This region coincides with the most stretched quadrilaterals. 

\begin{figure}
\includegraphics[width=1.\textwidth]{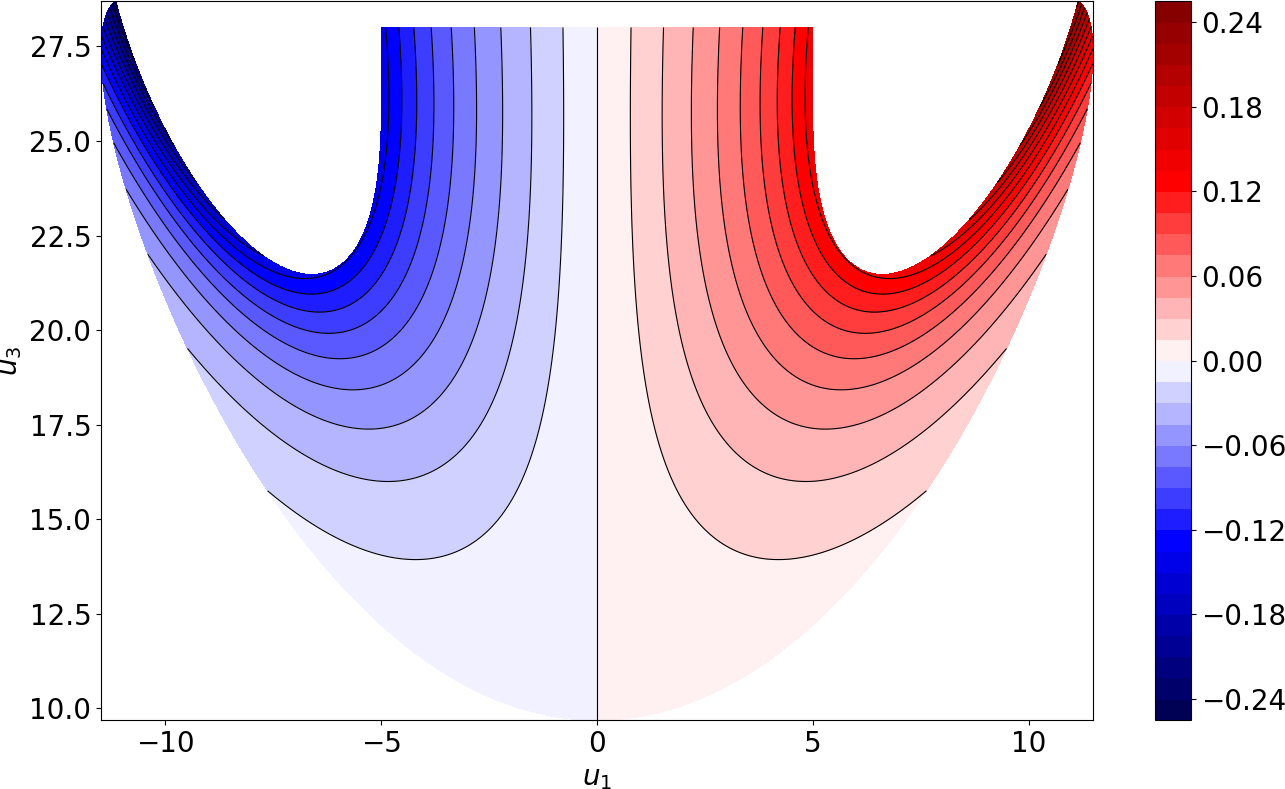}

\includegraphics[width=1.\textwidth]{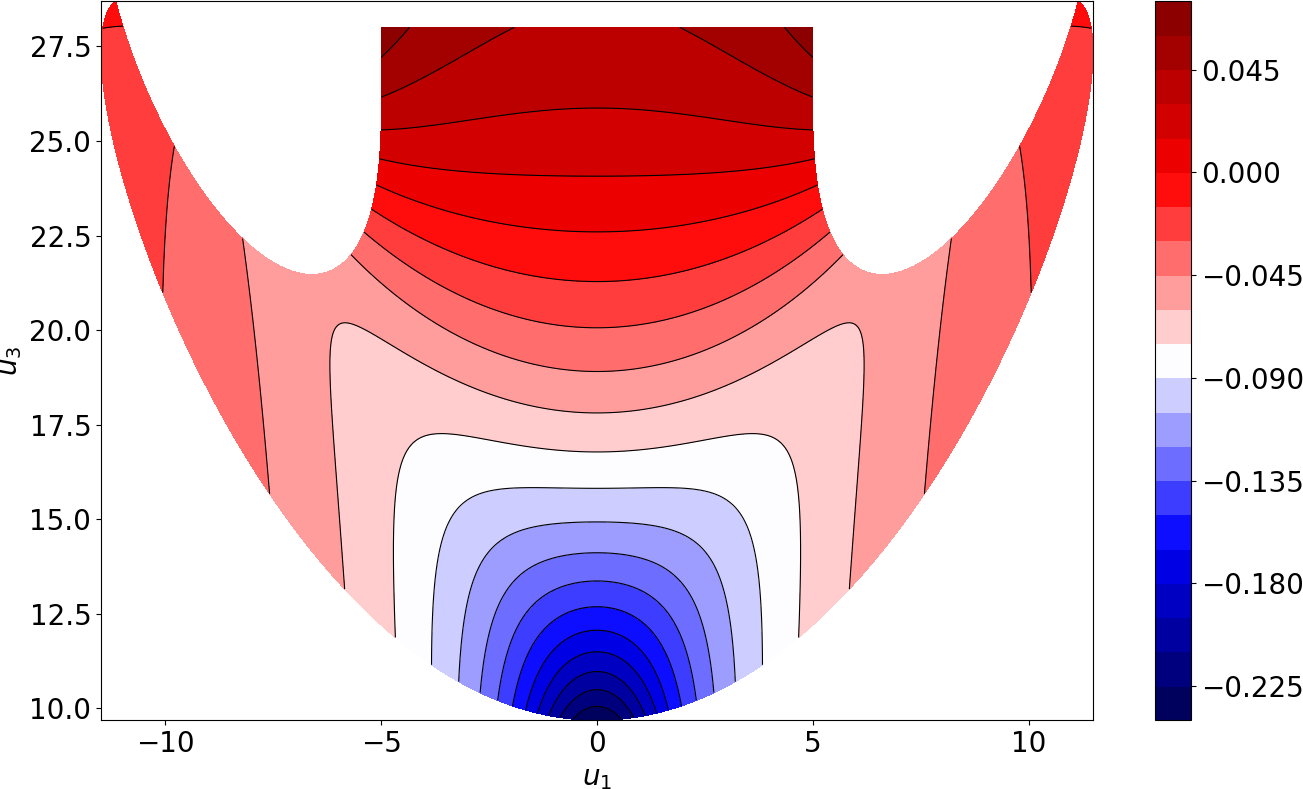}
\caption{$u_1\--u_3$ projection of the directional derivative of $\log\rho$, in the c-direction, $g_1:=g_c$ (upper plot), and $t$-direction, $g_2 := g_t$ (lower plot).}
\label{fig:general-density-gradient}
\end{figure}

\begin{figure}
\includegraphics[width=1.\textwidth]{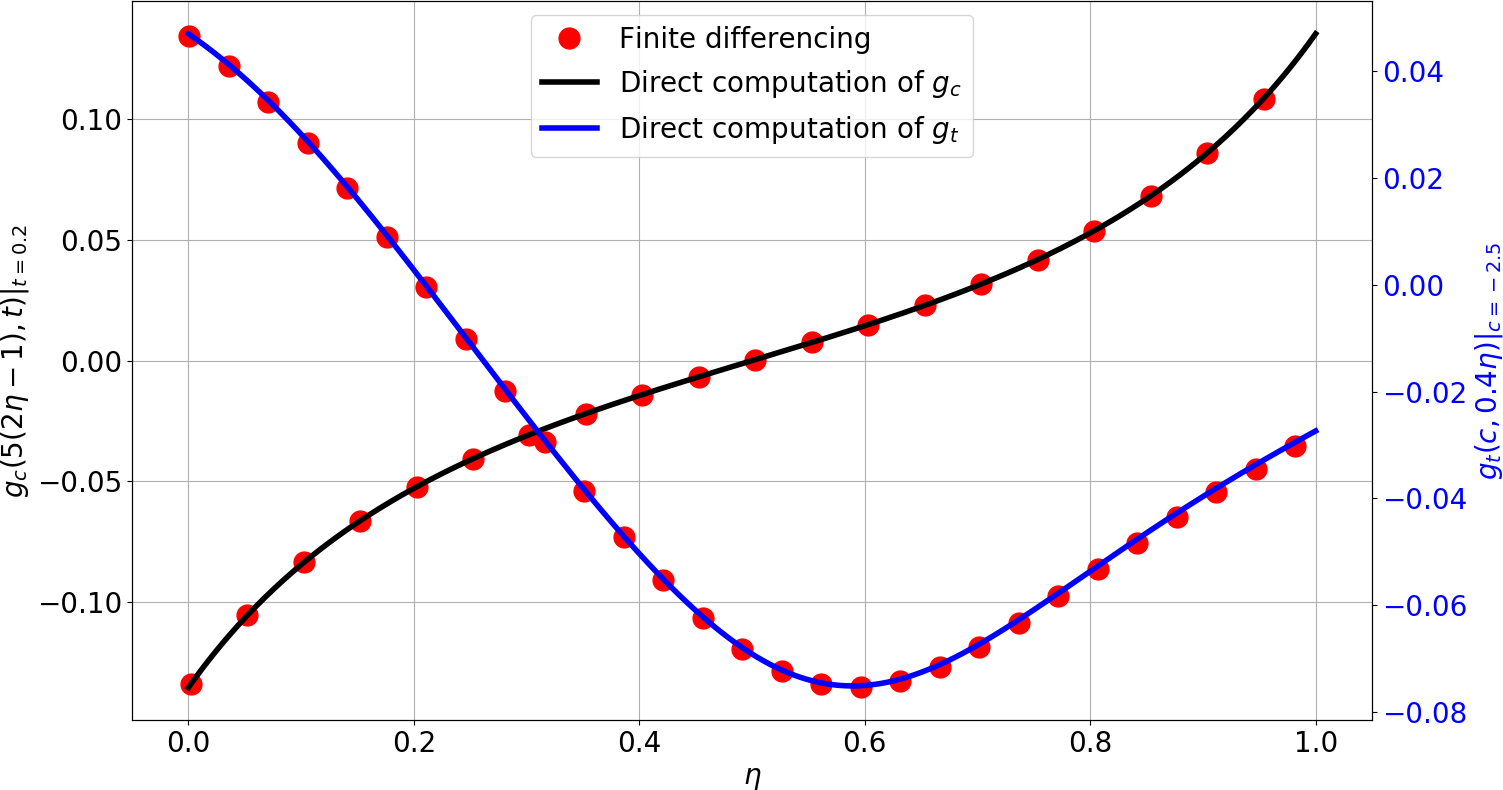}
\caption{The first and second component of the density gradient function $g=[g_c,g_t]$, respectively at $t=0.2$ (black solid line on Figure \ref{fig:general-mesh}) and  $c=-2.5$ (red solid line on Figure \ref{fig:general-mesh}), computed directly using Eq. \ref{eqn:general-g2-einstein} and through a finite difference method. In the latter approach, we note $\partial_{s_{i}} \rho = \partial_{\xi_{i}}\rho / \partial_{\xi_i} s_i$, $i = 1,2$, $\xi_1=c$, $\xi_2 = t$, where $s_i$ denotes the length of the isoparametric curve associated with $\xi_i$. Both the numerator and denominator is approximated using the central finite difference scheme on a uniform grid. The relation $s_{i}(\xi_i)$ is found in a way analogous to the one described in Section \ref{sec:curves}.}
\label{fig:general-direct-fd}
\end{figure}

Figure \ref{fig:general-density-gradient} shows the two components of the density gradient $g:=[g_c,g_t]^T = [\partial_{s_1}\log\rho,\partial_{s_2}\log\rho]^T $, corresponding respectively to the $c$- and $t$-direction. The distribution of $g_{c}$ is clearly symmetric with respect to the reflection points on the isoparametric line $c=0$, which is a manifestation of the fact the density is symmetric and directional derivative is computed in the direction of increasing $c$. The symmetry of $g_t$ is a direct consequence of the definition $g_t:=\partial_{s_2}\log\rho$, where $\log\rho$ itself is symmetric. Note the largest-in-magnitude values of $g_a$ concentrate around the boundaries of the range of $c$, i.e., at $c=\pm 5$ and, in case of $g_t$, around $u_1 = 0$. This reflects the fact the density gradient measures the relative rate of change of the density. In particular, its value becomes large if the rate of change of the density is large and/or the density itself is small. Figure \ref{fig:general-direct-fd} illustrates the density gradient along the bold isoparametric curves from Figure \ref{fig:general-mesh}, computed using Eq. \ref{eqn:general-g2-einstein} directly and through finite differencing. In case of both $g_c$ and $g_t$, we observe a good agreement between the solution computed directly and finite difference approximation, which validates our derivation of Eq. \ref{eqn:general-g2-einstein}.  

%SECTION 4------------------------
\section{Recursive algorithm for $g$ along trajectories defined by diffeomorphism $\varphi$}\label{sec:recursion}
Using the results presented in Section \ref{sec:1D} and \ref{sec:general}, we now propose an iterative method for the density gradient along trajectories defined by a $C^2$ diffeomorphism $\varphi:M^k\to M^{k+1}$, $k\in\mathbb{Z}$, where both $M^k$ and $M^{k+1}$ represent differentiable manifolds of the same dimension embedded in $\mathbb{R}^n$, $n\in \mathbb{Z}^{+}$. Let us consider two different charts, $x^{k}(\xi)\in N^k\subset M^k$ and $x^{k+1}(\xi)\in N^{k+1}\subset M^{k+1}$, such that
\begin{equation}
\label{eqn:recursion-map}
    x^{k+1}(\xi) = \varphi(x^k(\xi))
\end{equation}
for all $\xi\in V\subset{U}\subset\mathbb{R}^m$, $1\leq m\leq n$, $k\in\mathbb{Z}$. Let $\omega^k$ and $\omega^{k+1}$ be the natural volume forms in $M^k$ and $M^{k+1}$, respectively, where $\omega^{k+1}$ is the pushforward of $\omega^k$ under $\varphi$. Therefore, for all $k\in\mathbb{Z}$, the Lebesgue measure $m$ of the subspace $V$ can be expressed as follows,
\begin{equation}
    \label{eqn:recursion-lebesgue}
    m(V) = \int_{N^k}\rho^{k}(x)\;d\omega^k(x) =  \int_{N^{k+1}}\rho^{k+1}(x)\;d\omega^{k+1}(x),
\end{equation}
where $\rho^k$ and $\rho^{k+1}$ are densities implied by $x^k(\xi)$ and $x^{k+1}(\xi)$, respectively. Following the procedure involving Eq. \ref{eqn:general-measure}-\ref{eqn:general-metric}, it is possible to find the relation between $\rho^k$, $\rho^{k+1}$, and the metric tensors of the two transformations: $\xi\to x^k$ and  $\xi\to x^{k+1}$. Thus, by applying the chain rule, we find a relation between the parametric derivatives of $x^k(\xi)$ and $x^{k+1}(\xi)$, thanks to which a general recursive formula for the density gradient along the trajectory defined by $\varphi$ can be inferred. The $g^k$ function should be understood as the directional derivative of the (logarithmic) density implied by the chart $x^k(\xi)$.  In Section \ref{sec:recursion-derivation}, we derive an iterative procedure for $g^k$, while Section \ref{sec:recursion-example} presents the use of the proposed algorithm by revisiting the Lorenz '63 oscillator. Throughout this section, repeated indices in the subscript of any term imply summation (Einstein's notation), unless otherwise stated.    

\subsection{A generic recursive procedure for $g^k$}\label{sec:recursion-derivation}
As pointed out above, the first step is to find a relation between the parametric gradients of $x^k$ and $x^{k+1}$. Applying the definition of $\varphi$ from Eq. \ref{eqn:recursion-map} and the chain rule, we can expand $\nabla_{\xi}x^{k+1}$ in the following way,
\begin{equation}
    \label{eqn:recursion-first-deriv}
    \nabla_{\xi}x^{k+1}(\xi) = D\varphi(x^k(\xi))\;\nabla_{\xi}x^{k}(\xi),
\end{equation}
or, equivalently,
\begin{equation}
    \label{eqn:recursion-first-deriv-comp}
    \partial_{\xi_i}x^{k+1}(\xi) = D\varphi(x^k(\xi))\;\partial_{\xi_i}x^{k}(\xi),
\end{equation}
where $D\varphi$ denotes the $n \times n$ Jacobian matrix of $\varphi$, i.e., $(D\varphi)_{ij} = \partial_{x_j}\varphi_i$.
By differentiating Eq. \ref{eqn:recursion-first-deriv-comp} once more, with respect to $\xi_{j}$, we obtain
\begin{equation}
    \label{eqn:recursion-second-deriv-comp}
    \partial_{\xi_i}\partial_{\xi_j}x^{k+1}(\xi) = D^2\varphi(x^k(\xi))\left(\partial_{\xi_i}x^{k}(\xi),\partial_{\xi_j}x^{k}(\xi)\right) + D\varphi(x^k(\xi))\;\partial_{\xi_i}\partial_{\xi_j}x^{k}(\xi),
\end{equation}
where $D^2\varphi$ is the Hessian of $\varphi$, which is in fact a third-order $n\times n\times n$ tensor. Analogously to the example presented in Section \ref{sec:general-example}, the first term in the RHS of Eq. \ref{eqn:recursion-second-deriv-comp} is a bilinear form that outputs an $n$-element vector. In this case, the $i$-th component of that vector equals $\partial_{x_p}\partial_{x_q}\varphi_i(x^k(\xi))\;\partial_{\xi_i}x_p^k(\xi)\;\partial_{\xi_j}x_q^k(\xi)$. 

In the second step, we directly use the formula for the density gradient derived in Section \ref{sec:general-formula}. Let $f^k:=f(x^k(\xi))$ be a shorthand notation for any function $f$ defined at $x^k(\xi)$,  and $e_i(x^k(\xi)):= \partial_{\xi_i}x^k(\xi)$, $a_{ij}(x^k(\xi)) := \partial_{\xi_i}\partial_{\xi_j}x^{k}(\xi)$. Thus, by combining Eq. \ref{eqn:recursion-first-deriv-comp}, \ref{eqn:recursion-second-deriv-comp} with Eq. \ref{eqn:general-g2-einstein} derived for a generic chart $x(\xi)$, we conclude that 
\begin{equation}
    \label{eqn:recursion-g}
    g_i^k = - \frac{(R^{-1}_{lj})^k}{\|e_i^k\|}\,q_j^k\cdot a_{il}^k ,
\end{equation}
\begin{equation}
    \label{eqn:recursion-qr}
    (\nabla_\xi x)^k = [e_1^k\,e_2^k\dotsm e^k_m] = Q^k\,R^k = [q_1^k\,q_2^k\dotsm q^k_m]\,R^k.
\end{equation}
\begin{equation}
    \label{eqn:recursion-e}
    e^{k+1}_{i} = D\varphi^k e^k_i,
\end{equation}
\begin{equation}
    \label{eqn:recursion-a}
    a^{k+1}_{ij} = D^2\varphi^k (e_i^k,e_j^k) + D\varphi^k a^k_{ij},
\end{equation}
hold for any $\xi\in{V}\subset U$. 

To summarize, if a map $\varphi$ relating two consecutive points on the trajectory, $x_{k}(\xi)$ and $x_{k+1}(\xi)$, is available, then the density gradient at one point can be computed using information associated with the other point. In particular, according to Eq. \ref{eqn:recursion-g}, the $i$-th component of $g$ requires knowledge of $e_j$, $j = 1,...,m$ and $a_{pq}$,  $p, q = 1,...,m$ at the same point. Thus, to compute one component of the density gradient at $x^k(\xi)$ for some $\xi$, we need to apply the recursion in Eq. \ref{eqn:recursion-e} $km$ times and, analogously, the recursion in Eq. \ref{eqn:recursion-a} $1/2\,k m^2$ times. The $1/2$ factor is a consequence of the fact that $a_{ij}(\xi) = a_{ji}(\xi)$ for any admissible $\xi$, because $x_{k}$ is assumed to be twice differentiable for any $k\in\mathbb{Z}$. In addition, at every step $k$, the QR factorization of $(\nabla_{\xi} x)^k = [e_1^k\,e_2^k\dotsm e^k_m]$ and inversion (either direct if $m$ is small or through solving a linear system) of the resulting $m\times m$ $R^k$ matrix must be performed. We assume $x^0(\xi)$ is given, from which we directly compute initial conditions for recursions in Eq. \ref{eqn:recursion-e} and \ref{eqn:recursion-a}.

The recursion involving Eq. \ref{eqn:recursion-g}-\ref{eqn:recursion-a} can be used to devise algorithms for differentiating the invariant, physical SRB measure $m_{\mathrm{SRB}}$, which is guaranteed to exist in uniformly hyperbolic systems. In general, $m_{\mathrm{SRB}}$ is not absolutely continuous everywhere on the manifold, but only conditional measures of $m_{\mathrm{SRB}}$ along unstable manifolds are absolutely continuous. The SRB density gradient $g_{\mathrm{SRB}}$, defined as a directional derivative of the conditional SRB density on the unstable manifold, is a byproduct of the integration by parts (analogous to Eq. \ref{eqn:intro2}), preceded by the disintegration of $m_{\mathrm{SRB}}$ \cite{chandramoorthy-s3,sliwiak-1d}. Thus, if a direction of the unstable manifold is given, the recursive formula presented in this section might be further developed to compute the SRB density gradient, defined on a manifold of any dimension, along a trajectory initiated at a $m_{\mathrm{SRB}}$-typical point. 

\subsection{Example: evolution of a 1D manifold}\label{sec:recursion-example}

In this section, we demonstrate the application of the recursive scheme for the density gradient $g^k$. For this purpose, let us re-consider the Lorenz '63 oscillator, defined by System \ref{eqn:general-lorenz}. In particular, we define $\varphi$, such that it represents numerical time integration of System \ref{eqn:general-lorenz} for a period of $\Delta t$, i.e., $u(t+\Delta t) = \varphi(u(t))$ with $u(t)$ being the solution of the system at time $t$. Let us consider a 1D smooth manifold embedded in $\mathbb{R}^3$ described by the following chart $x^0(c) = [c,c,28]^T$, $-5\leq c\leq 5$. Note $x^0(c)$ coincides with the black solid boundary of the surface depicted in Figures \ref{fig:general-mesh}-\ref{fig:general-mesh-ext}. Now, by applying $\varphi$ recursively, the next step is to numerically compute a sequence of charts $\{x^0(c), x^1(c), x^2(c), ...\}$, where $x^{k+1}(c) = \varphi(x^{k}(c))$. Our aim is to compute $g^k = \partial_s\log\rho^k$, where $\rho^k$ is a density implied by the chart $x^k(c)$. The operator $\partial_s$ denotes a generic directional derivative along the curve in the direction of increasing $c$. The formulas derived in Section \ref{sec:recursion-derivation} give us all necessary tools to compute $g^k$ along the trajectory defined by $\varphi$. In this example, however, we consider the simplest case, $m=1$. Eq. \ref{eqn:recursion-g}-\ref{eqn:recursion-a} can be dramatically simplified, because $\nabla_{\xi} x = dx/dc$ is just a vector, and thus QR factorization is equivalent to normalizing that vector. Let $e = dx/dc = \|dx/dc\|\,q$ and $a = d^2x/dc^2$ and, therefore,
\begin{equation}
    \label{eqn:recursion-g-1D}
    g^k = - \frac{q^k\cdot a^k}{\|e^k\|^2},
\end{equation}
\begin{equation}
    \label{eqn:recursion-e-1D}
    e^{k+1} = D\varphi^k\,e^k,\;\;\;q^k = \frac{e^k}{\|e^k\|}, 
\end{equation}
\begin{equation}
    \label{eqn:recursion-a-1D}
    a^{k+1} = D^2\varphi^k(e^k,e^k) + D\varphi^k\,a^k. 
\end{equation}
Note Eq. \ref{eqn:recursion-g-1D}-\ref{eqn:recursion-a-1D} can be derived directly using Eq. \ref{eqn:curve3} and the chain rule for parametric derivatives.

\begin{figure}
    \centering
    \includegraphics[width=1.\textwidth]{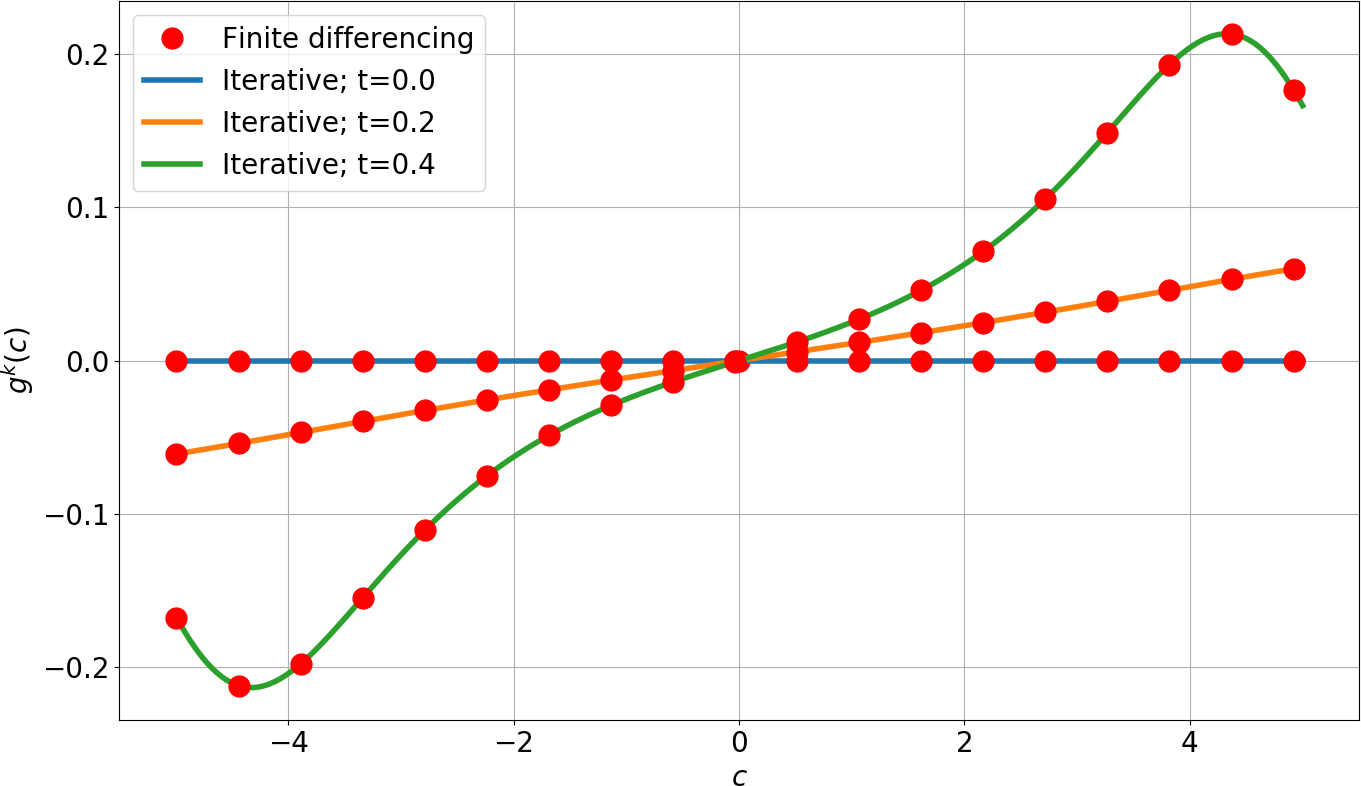}    
    \caption{Density gradient function computed using the recursion involving Eq. \ref{eqn:recursion-g-1D}-\ref{eqn:recursion-a-1D} at three different time steps $k = t/\Delta t$. The finite difference approximation is generated using the approach described in Section \ref{sec:curves}.}
    \label{fig:recursion-fd}
\end{figure}

How does this example differ from the one presented in Section \ref{sec:general-example}? There, we used a chart $x_s(\xi):\mathbb{R}^2\to\mathbb{R}^3$, $\xi=[c,t]^T$, which defined a two-dimensional manifold. The rate of change of $x(\xi)$ in the $t$-direction was determined by the Lorenz '63 oscillator (System \ref{eqn:general-lorenz}). Here, using the iterative procedure, we generate a bunch of 1D manifolds $x^k(c)$. The evolution of these curves (in geometric sense) is determined by $\varphi$, which is in fact a discrete version of System \ref{eqn:general-lorenz}. Thus, if we generate infinitely many such curves and $\Delta t\to 0$, we effectively obtain the same surface as the one shown in Figure \ref{fig:general-mesh}. Intuitively, the density $\rho_s$ implied by $x_s(\xi)$ measures number of points mapped from a uniform distribution per unit surface area. Likewise, the density $\rho^k$, implied by the chart $x^k(c)$, measures number of points mapped from a uniform distribution per unit curve length. Since we use the same discretization scheme to integrate differential equations, the localization of points obtained in both the computation of surface from Section \ref{sec:general-example} and, here, evolution of curves is exactly the same. However, the density $\rho^k$ does not equal to the marginal distribution of $\rho_s$ at $t = k\,\Delta t$ (assuming uniform discretization of time). In case of the surface example, the value of the density function reflects the densification of points, mapped from a uniform distribution, in both the $t$ and $c$ directions. In the latter example, the density is determined only by the localization of points along the evolving curve. Figure \ref{fig:recursion-fd} illustrates the density gradient $g^k$ along the evolving curve, recorded at three different time steps $k$. We observe $g^k=0$ at $k=0$, which is a consequence of the choice of the uniformly distributed initial condition. Due to the symmetric geometry of $M^k$, defined by the Lorenz '63 oscillator at $t\in[0,0.4]$, the density gradient features symmetric behavior with respect to the origin of the $g^k(c)$-vs.-$c$ relation. 

\section{Conclusions and future work}\label{sec:conclusion}
A Monte Carlo integration scheme applied to a highly-oscillatory function might be remarkably expensive. The computational cost, however, can be dramatically reduced by integrating the original formulation by parts. Such treatment gives rise to a new quantity, i.e., a directional derivative of the logarithm of the density implied by a chart describing the integration domain. The computation of that derivative, which we call the {\em density gradient}, requires knowledge of the first and second derivatives of the chart with respect to the domain parameterization. If the domain manifold evolves according to some diffeomorphism $\varphi$, the calculation of the density gradient along a trajectory requires solving a collection of first- and second-order tangent equations, involving both the Jacobian and Hessian of $\varphi$. The number of these equations is respectively proportional to $m$ and $m^2$, where $m$ is the dimension of the manifold. 

The formulas derived in this paper is a major step toward constructing generalizable algorithms for {\em SRB density gradient}. This quantity plays a major role in the sensitivity analysis of uniformly hyperbolic systems, including many popular chaotic equations. Using the recursive formula for the density gradient along a trajectory defined by $\varphi$ and the definition of the SRB measure, one can potentially devise a trajectory-driven procedure for the SRB density gradient. This is in fact the subject of the authors' ongoing investigation.  

\section*{Acknowledgments}
This work was supported by Air Force Office of Scientific Research Grant No. FA8650-19-C-2207 and U.S. Department of Energy Grant No. DE-FOA-0002068-0018.

\bibliographystyle{elsarticle-num-names}
\bibliography{references.bib}

\end{document}